\documentclass[12pt,reqno,a4paper]{amsart}
\usepackage{amsmath,amssymb,amsfonts,amscd}
\usepackage[mathscr]{eucal}

\usepackage{comment}
\usepackage{color}
\usepackage{hyperref,graphicx}
\usepackage[normalem]{ulem} 
\usepackage{tikz}
\usetikzlibrary{calc, arrows.meta}
\usepackage{arydshln}
\usepackage{tikz-cd} 

\allowdisplaybreaks

\date{}
 \thanks{This work was partially supported by Simons Foundation collaboration grant for mathematicians number 846970}    

\title[Berezinian expansion and super exterior powers]{Berezinian expansion and super exterior powers}
\author{Maheshan~Ekanayaka and Ekaterina~Shemyakova}
\address{Department of Mathematics,  University of Toledo, Toledo,  Ohio, USA}
\email{ekaterina.shemyakova@utoledo.edu and mekanaya@gettysburg.edu}

\newtheorem{theorem}{Theorem}[section]

\newtheorem{lemma}{Lemma}[section]
\newtheorem{corollary}{Corollary}[section]
\theoremstyle{definition}

\newtheorem{example}{Example}[section]
\newtheorem{remark}{Remark}[section]
\newcommand{\R}{\mathbb{R}}

\newcommand{\Z}{\mathbb{Z}}

\newcommand{\ZZ}{{\mathbb Z}}
\DeclareMathOperator{\diag}{diag}
\DeclareMathOperator{\Ber}{\mathop{Ber}}
\DeclareMathOperator{\Vol}{\mathop{Vol}}
\DeclareMathOperator{\tr}{\mathop{Tr}}
\DeclareMathOperator{\str}{\mathop{Str}}

\providecommand{\C}{\mathbb{C}}

\newcommand{\e}{{\varepsilon}}

\theoremstyle{definition}

\def\Xint#1{\mathchoice
{\XXint\displaystyle\textstyle{#1}}%
{\XXint\textstyle\scriptstyle{#1}}%
{\XXint\scriptstyle\scriptscriptstyle{#1}}%
{\XXint\scriptscriptstyle\scriptscriptstyle{#1}}%
\!\int}
\def\XXint#1#2#3{{\setbox0=\hbox{$#1{#2#3}{\int}$ }
\vcenter{\hbox{$#2#3$ }}\kern-.6\wd0}}

\def\dashint{\Xint-}
\newcommand{\fint}{\dashint}

\begin{document}
\begin{abstract}In supergeometry, the classical identification of top-degree differential forms with volume forms suitable for integration breaks down. To address this issue, several generalized notions of differential forms have been introduced, including integral forms and pseudo-differential forms of Bernstein-Leites and $r|s$-forms of Voronov-Zorich. The Baranov-Schwarz transformation maps pseudo-differential forms to $r|s$-forms of type II, and, on a supermanifold of dimension $n|m$, integral $r$-forms are isomorphic to $r|m$-forms. However, an explicit construction of $r|s$-forms for arbitrary $s$ has remained elusive.
	
First, we observe that this problem can be related with expansions of $\Ber(E+zA)$ for a linear operator $A$ on a superspace $V$. It was already known that the coefficients of the expansions near zero and near infinity give the supertraces of the representations $\Lambda^{r|s}(A)$ in the two extreme cases $s=0$ and $s=m$, respectively. 

In this paper we show that the intermediate expansions of $\Ber(E+zA)$, taken in the annular regions between consecutive poles, encode the supertraces of representations on explicitly constructed vector spaces that model $\Lambda^{r|s}(V)$ for $0<s<m$. We introduce a formal analogue of the Baranov-Schwarz transformation and show that it relates these spaces to Voronov-Zorich $r|s$-forms.
	
As a separate remark, we establish  that $1|1$-forms  at a point of a supermanifold of dimension $n|m$ can be realized as closed differential forms on the super projective space $\mathbb{P}^{m-1|n}$. 
	
\end{abstract}

\maketitle

%\tableofcontents

\section{Introduction}

In the purely even case (for ordinary manifolds), top-degree differential forms serve for integration on manifolds as volume forms. The direct analogue of this identification breaks down in the supergeometric setting. 

Indeed, on a supermanifold $M = M^{n|m}$, one can define the de Rham complex consisting of differential forms locally generated by the differentials
$dx^i$ and $dx^{\hat{j}}$ of the even and odd coordinates, respectively. However, this direct analog of the classical construction presents several issues.
Here $d$ is an odd operator, and the differentials $dx^{\hat{j}}$ of the odd coordinates are themselves even, which means there are arbitrarily high-degree monomials in $dx^{\hat{j}}$. Moreover, under a change of coordinates, the number of $dx^i$ and the number of $dx^{\hat{j}}$ are not separately preserved, so a bi-grading is not well-defined. Only the total number of all $dx^i$ and $dx^{\hat{j}}$ combined is preserved and is used for grading. The third issue is that differential forms constructed in this way do not transform as densities under  changes of coordinates. All of that makes them unsuitable for Berezin integration.

Among the key ideas addressing this problem are:
\begin{enumerate}
	\item \emph{Pseudo-differential forms}, introduced by Bernstein and Leites~\cite{Bern_Lei77_part2}, are defined as elements of $C^\infty (\Pi TM)$. They are often considered under the additional condition of being compactly supported or   rapidly decreasing as functions of $d x^{\hat j}$. These can be integrated over supermanifolds, but they do not have  natural grading.
	\item \emph{Integral forms}, introduced by Bernstein and Leites~\cite{Bern_Lei77_part1}, are fiberwise polynomial functions on $\Pi T^*M$ whose coefficients are Berezin volume forms. A typical element has the form  
	\begin{equation*}
		\sigma = f(x) \, D(x) \cdot (x^*_1)^{\alpha_1} \cdots (x^*_n)^{\alpha_n} (x^*_{\hat{1}})^{\beta_1} \cdots (x^*_{\hat{m}})^{\beta_m} \, ,
	\end{equation*} 
	where $x^*_i$ and $x^*_{\hat j}$ denote the fiber coordinates on $\Pi T^*M$.
	Integral forms of degree $r$ are integrable over 
	subsupermanifolds of dimension $r|m$, where $m$ is the maximal odd dimension. Here $r=n-\alpha_1-\ldots-\alpha_n-\beta_1-\ldots-\beta_m$.
	\item \emph{$r|s$-forms}, introduced and developed by Voronov and Zorich in~\cite{tv:compl,tv:pdf,tv:cohom-1988,tv:dual}, 
	are defined as $r|s$-densities $L(x, \dot{x})$ that satisfy a system of partial differential equations~\eqref{eq.fund}.	An $r|s$-form can be integrated over a ``surface'' or ``singular supermanifold'' of dimension $r|s$,  $0\leq s \leq m$. There is a remarkable relation with integral geometry in the sense of Gelfand-Gindikin-Graev~\cite{gelfand-gindikin-graev:intgeometry1980}. (See~\cite{tv:pdf,tv:cohom-1988},  \cite{gel:gms}.)
\end{enumerate}
The \emph{Baranov–Schwarz transformation}~\cite{Baranov_Shwarz} maps pseudo-differential forms into $r|s$-forms of type II (see~\cite{tv:pdf}). Also, integral $r$-forms are isomorphic to $r|m$-forms. We recall these notions  in greater detail  in sections~\ref{Sec:VZ}, \ref{sec:BerLei} and \ref{Sec:BSch} below.

It is an open problem to construct $r|s$-forms  or the corresponding    $r|s$-exterior powers $\Lambda^{r|s}(V)$  for arbitrary $s$ explicitly. (Meaning general description, not particular examples.)

In the present paper, we obtain two results in this direction.

First, we show that on a supermanifold of dimension $n|m$, the vector space of $1|1$-forms at a point is isomorphic to the space of closed  differential forms on the super projective space $\mathbb{P}^{m-1|n}$ (which is the projectivization of $\Pi T_xM$). %Algebraically, this means the description of $\Lambda^{1|1}(V)$. 

Secondly, we  consider the following related problem. % involving $r|s$-forms  for arbitrary $r$ and $s$. 

For an  even operator $A$ on a vector superspace $V$ of dimension $n|m$,   
it was known that the function $\Ber(E + z A)$ serves 
as a kind of generating function for the supertraces $\str\Lambda^r(A)$ of the representations  in the exterior powers $\Lambda^r(V)$, and, after certain manipulation, also for the supertraces  $\str\Lambda^{r|m}(A)$ in the spaces $\Lambda^{r|m}(V)$ similar to integral forms; see~\cite{Kh_V_recurrent}. In more detail, $\Ber(E + z A)$ generates these supertraces via its power expansions in $z$ at zero and at infinity.
Notice that the spaces $\Lambda^{r}(V)$ and $\Lambda^{r|m}(V)$ (where the odd degree is maximal) represent extreme cases. Thus, the expansion of $\Ber(E + zA)$ is related to these extreme cases of exterior powers. Since intermediate degrees also occur, a natural question arises.

Indeed, $\Ber(E + z A)$ as a rational function also admits expansions in the annular regions between its poles. A natural --- and solved in this paper --- question is whether these expansions can be connected to the supertraces  $\str\Lambda^{r|s}(A)$ %for $s \neq 0, m$. By $\Lambda^{r|s}(A)$ we mean the representation 
of the representations in the spaces $\Lambda^{r|s}(V)$, for $s \neq 0, m$. One is tempted to think that there must be a connection because the number of possible types of such spaces $0\leq s\leq m$ miraculously coincides with the number of the annular regions including neighborhoods of zero and infinity, which correspond to $s=0$ and $s=m$.  

One of the challenges lies in the fact that, for fixed $r$, in the previously studied cases the spaces $\Lambda^{r|0}(V)$ and $\Lambda^{r|m}(V)$ are finite-dimensional, whereas the spaces $\Lambda^{r|s}(V)$ for $s \neq 0, m$ have to be infinite-dimensional. Thus, there is a problem of properly defining and understanding  supertraces in this infinite-dimensional context; particularly, in the absence of an explicit description of these spaces. 

In the present paper, we show that the coefficients of the expansions of $\Ber(E + z A)$ in the  intermediate  annular regions 
are equal to the supertraces on certain vector spaces that we explicitly construct and which, informally speaking, are candidates for the super exterior powers $\Lambda^{r|s}(V)$. (For $s=0$ and $s=m$ the same construction gives spaces isomorphic to the known $\Lambda^{r|0}(V)$ and $\Lambda^{r|m}(V)$.) The objects from which these spaces are constructed  --- formal delta-functions and their derivatives --- arose before in string theory (Belopolsky~\cite{belopolsky1997new} and Witten~\cite{Witten_Notes_on_Super}). Under a formal Baranov-Schwarz transformation that we introduce, they correspond to ``Pl{\"u}cker forms'' or ``non-linear wedge products'' considered in~\cite{shemya:voronov2019:superplucker}.  There are questions that remain, such as dependence of a choice of basis, and that will be addressed elsewhere. We also  plan to develop a version of super de Rham theory based on our construction, which will   incorporate $r|s$-forms with $r<0$~\cite{tv:dual}.

The structure of the paper is as follows. 

In Section~\ref{Sec:E}, we recall the known facts about the expansions of $\Ber(E + z A)$ at zero and infinity and show an example of the expansion in the intermediate region, which will be the central topic in Sections~\ref{sec.more} and~\ref{sec.geom}.  In Section~\ref{Sec:VZ}, we recall Voronov-Zorich $r|s$-forms, paying particular attention to the cases where an explicit description   exists. In Sections~\ref{sec:BerLei} and~\ref{Sec:BSch}, we recall Bernstein-Leites pseudo-differential forms and integrals forms, and the Baranov-Schwarz integral transformation.

In Section~\ref{Sec:11},  
we consider the space of $1|1$-forms at a point (or, algebraically, $\Lambda^{1|1}(V)$) and show that it is isomorphic to the subspace of all closed $1$-forms on $\mathbb{P}^{m-1|n}$ (Theorems~\ref{thm.11} and \ref{thm.11alg}). This is our first main result. There is no explicit relation with integral geometry here, but it is tempting to think that it exists.

Sections~\ref{sec.more} and~\ref{sec.geom} are central in this paper. In Sections~\ref{sec.more}, we calculate the Laurent expansions of the function $\Ber(E + z A)$ in all regions in terms of the eigenvalues of $A$. %  (in particular, near zero and near infinity, thus re-proving the theorems of Khudaverdian-Voronov and Th.~Schmitt). 
Geometric interpretation of these formulas  follows  in Section~\ref{sec.geom}. We consider the parity-reversed  space $\Pi V$ and define vector spaces   $S^{N,s}(\Pi V)$.  
(Note ``generalized  symmetric powers'' instead of  exterior powers.) 
We prove that in the region between  consecutive poles (including the cases of near $0$ and near $\infty$), the coefficient of the expansion of $\Ber(E + z A)$   at $z^N$ is the supertrace of the representation of $A$ in   $\Pi^NS^{N,s}(\Pi V)$ (Theorem~\ref{thm.geomint}).  This is our second main result. Here   $-\infty<N<+\infty$ for  $0< s< m$   (the intermediate case, not known before), and $N\geq 0$ or $N\leq n-m$  for $s=0$ and $s=m$, respectively. % There we explicitly, in terms of bases, construct spaces that are ``candidates for super exterior powers'' for all $0\leq s\leq m$ in dimension $n|m$. 

In Section~\ref{sec.conVZ}, we introduce   \emph{formal Baranov-Schwarz transformation} and analyze connection of the spaces $S^{N,s}(\Pi V^*)$ with   Voronov-Zorich forms (where $N=r-s$).

\emph{Acknowledgement.} 
We thank Th. Voronov for very useful discussions. 
We thank A. Odesskii for suggesting an  interpretation of formal delta-functions in terms of power series (see in~\ref{sec.fdeltas}).

\emph{Notation.} We use standard supergeometric notation (see e.g. Manin~\cite{manin:gaugeeng}). In particular, $\Pi$ denotes the parity reversion functor for modules and vector bundles. 
%
%1. Consider interpretation $\omega_L$ interpretation for $1|1$ in 
%$n|m$. Same about $r|1$. Will they be closed?
%
%2. Consider Laurent polynomials of my type. How can they be expanded? 
%What will happen if we do a change of variables in the ambient space?
%
%3. Investigate rational functions with fixed poles.
% 
%4. Consider 
%
%%However, this is not the case when one naively generalizes the standard notion of differential forms to the supermanifold setting. 
%
%In the pure even case, differential forms have a well-defined degree, with only top-degree forms (with the degree equal to the manifold's dimension) can be integrated over the entire manifold. These top forms transform like volume elements (i.e., densities or measures) under orientation-preserving coordinate changes -- so integration makes geometric sense.
%
%
%In the supermanifold case, where coordinates are split into even coordinates $x^a$ and odd coordinates $\xi^i$,  a naive generalization leads to forms with odd symbols $dx^a$,
%and even symbols $d\xi^i$,
%where the $d\xi^i$ commute without a sign, allowing forms of arbitrary (even infinite) degree. These generalized differential forms \textbf{do not} transform like Berezinian volume elements. As a result, they are \textbf{not suitable for integration} (see details in Subsec.~\ref{Subsec:BerLei}).
%
%Belopolsky~\cite{belopolsky1997new}

\section{Expansions of the Berezinian}
\label{Sec:E}

Consider an even linear operator $A$ in an $n|m$-dimensional superspace $V$\,\footnote{Strictly speaking, we should understand by $V$ an $n|m$-dimensional free module over some auxiliary Grassmann algebra.}. We will consider the power expansions of the rational function $\Ber(E + z A)$. Following~\cite{Kh_V_recurrent}, we refer to it as the \emph{characteristic function} of $A$. 

In the   purely even case, $\Ber$ is $\det$  and we have the  polynomial   $\det(E+zA)$. It  differs from the characteristic polynomial $\det(A-zE)$   only by renumbering    terms and   signs. (In the supercase,   $\Ber(E + z A)$   is more convenient  than   $\Ber(A-zE)$.)
The expansion formula for $\det(E+zA)$,
\begin{equation} \label{eq:detext}
\det (E + z A)=1 + z\tr \Lambda^1(A)  + z^2\tr\Lambda^2(A)  + \ldots + 
z^{n-1}\tr\Lambda^{n-1}(A)  +
z^n\det(A)       
\end{equation} 
has appeared in   advanced linear algebra  books since the first decades of the 20th century; see, for example, Weyl~\cite{weylbook}. The first explicit mentioning in terms of the traces of the exterior powers may be attributed to E. Cartan. % and became standard after Bourbaki.

In 1981, Schmitt~\cite{Schmidt_recurrence} showed that the polynomial expansion~\eqref{eq:detext} for   determinant in the purely even case can be extended to an analogous expansion for   Berezinian:
 \begin{equation} \label{exp_zero}
 	\Ber (E + z A)= 	\sum_{k = 0}^{\infty} c_{k} z^{k}  = 1 + \str \Lambda^1(A)  z + \str \Lambda^2(A)   z^2+ \dots \, ,
 \end{equation}
 where $\str$ denotes the supertrace.
Since transitioning from determinants in the purely even case to Berezinians in the super case introduces rational functions instead of polynomials, it is useful to note that expansion~\eqref{exp_zero} is a Taylor series expansion at zero.
 
Later, in 2005, Khudaverdian and Voronov~\cite{Kh_V_recurrent} 
considered the expansion at infinity and found the geometric meaning for the coefficients:  
\begin{equation} \label{exp_infty}
	\mathrm{Ber}(E + zA) = 	\sum_{k = -\infty}^{n-m} c_{k}^{*} z^{k} \quad \text{where } c_{k}^{*} = \mathrm{Ber}\, A \cdot \str \, \Lambda^{n - m - k} (A^{-1}) 
	=\str \Sigma^{k+m} (A) \, .
\end{equation}
Here $\Sigma^{k+m} (A)$ denotes the induced action of $A$ in the vector space $\Sigma^{k+m}(V):=\Ber V\otimes \Lambda^{n-m-k}(V^*)$, $\Sigma^{k+m} (A)=\Ber A\cdot \Lambda^{n-m-k}({A^{-1}}^*)$.

The non-finiteness of the above series is compensated by its coefficients satisfying a recurrence relation from a certain index onward.
    
\begin{theorem}[\cite{Kh_V_recurrent}] \label{thm:REC-REL_ZERO} 
Let $\dim V=n|m$. 
\begin{enumerate}
	\item The coefficients $c_{k}(A) = \str \Lambda^{k}(A)$ of the expansion of $\Ber(E+zA)$ at zero~\eqref{exp_zero} satisfy a recurrence relation of length $m+1$:
\begin{equation}
c_k b_m+c_{k+1}b_{m-1}+ \dots +c_{k+m-1}b_{1}+c_{k+m}b_{0}=0 \ ,	
\end{equation}
for all $k > n-m$.
 Here $b_0=1$, $b_1, \dots, b_{m}$ are the same 
in every relation (do not depend on $k$). 
Here, in the case $n < m$, we set $c_{k}=0$ for $n-m < k < 0$.
%That the coefficients in the expansion of a rational function satisfy a linear recurrence relation is classical, and this property characterizes rational functions (Kronecker’s theorem). The significance of Theorem~\ref{thm:REC-REL_ZERO} is that such a recurrence is satisfied by the supertraces of the exterior powers -- serving as a substitute for the vanishing of higher exterior powers in the ordinary (non-super) case.
\item
The coefficients $c_{k}^{*}(A) = \str \Sigma^{m+k}(A)$ of the expansion of $\Ber(E+zA)$ at infinity~\eqref{exp_infty} satisfy the same recurrence relation:
\[
 c_{k-m}^{*} b_{m} + \dots + c_{k}^{*} b_{0} = 0 
\]
for all $k < 0$.   Here, in the case $n<m$, we set $c_{k}^*=0$ for $n-m < k < 0$.

\item 
Sequences $c_{k}$ and $c_{k}^{*}$ can be combined into a single recurrent sequence by considering the differences:
\[
\gamma_{k} = c_{k} - c_{k}^{*},
\]
which satisfy the recurrence relation:
\[
b_{0}\gamma_{k+q} + \cdots + b_{q}\gamma_{k} = 0 
\]
for \emph{all} values of $k \in \Z$.
\end{enumerate}
\end{theorem}

\begin{example} Let $\dim V =1|2$ and fix a homogeneous basis in which a linear operator 
	$A$ is represented by a diagonal matrix:
\begin{equation}
	A = \left(
\begin{array}{c|cc}
	x & 0 & 0 \\
	\hline
	0 & y_1 & 0 \\
	0 & 0 & y_2
\end{array}
\right)_{\!\!} \, ,
\end{equation}
where $x$, $y_1$, $y_2$ are even. Assume $|y_2|<|y_1|$. Then, we obtain:
\begin{equation}
	\Ber(1+zA) =
\frac{1+zx}{(1+zy_1)(1+zy_2)}=\frac{A}{1+zy_1}+\frac{B}{1+zy_2}  \, , 
\end{equation}
where $A=\frac{x-y_1}{y_2-y_1}$ and $B=\frac{y_2-x}{y_2-y_1}$.
    
\begin{figure}[h]
\begin{tikzpicture}[scale=2, >=Stealth]

% Axes
\draw[->] (-2, 0) -- (2, 0) node[right] {$\mathrm{Re}(z)$};
\draw[->] (0, -1.5) -- (0, 1.5) node[above] {$\mathrm{Im}(z)$};

% Poles
\node[below left] at (0,0) {$0$};
% \filldraw[red] (0,0) circle (1pt);

% Poles at z = -1/y₁ and z = -1/y₂
\def\yoneval{2} % y₁ = 2 (pole at z = -1/2)
\def\ytwoval{1} % y₂ = 1 (pole at z = -1)

\filldraw[black] (-1/\yoneval, 0) circle (1pt) node[below left] {$-\frac{1}{|y_1|}$};
\filldraw[black] (-1/\ytwoval, 0) circle (1pt) node[below left] {$-\frac{1}{|y_2|}$};

% Convergence regions (circles with positive radii)
\draw[dashed, black] (0,0) circle (1/\yoneval); % |z| < 1/y₁
\draw[dashed, black] (0,0) circle (1/\ytwoval); % |z| < 1/y₂

% Labels
% \node[green, above right] at (0.5, 0.1) {$|z| < \frac{1}{y_1}$};
% \node[orange, below right] at (1, -0.1) {$|z| < \frac{1}{y_2}$};

\end{tikzpicture}

\end{figure}

Expansion at zero: in the region $|z|<\frac{1}{|y_1|}$, we have 
%the following expansion:
\begin{align}
\frac{1+zx}{(1+zy_1)(1+zy_2)}
%=A+B-z(A y_1+B y_2)+z^2(A(y_1)^2+B(y_2)^2)-z^3(A(y_1)^3+B(y_2)^3)+\dots \\    
&
%=\sum_{k=o}^{\infty}(z)^{k}(-1)^{k}(A(y_1)^k+B(y_2)^k)
=\sum_{k=0}^{\infty}z^{k}c_{k} \, , \\
c_k &= (-1)^{k}(A(y_1)^k+B(y_2)^k) \, , \ k \geq 0 \, , \label{eq:ck_closed_formula}
 \end{align}
 where  $c_k$ satisfy the following recurrence relation:
\begin{align}
      &  c_k y_1y_2 + c_{k+1} (y_1+y_2) + c_{k+2} = 0 \, , \ k \geq 0 \, ,  \\
      & c_0 = 1 \, , \ c_1 =x-y_1 -y_2 \, .
\end{align}
Expansion at infinity: in the region $\displaystyle |z|> \frac{1}{|y_2|}$, we have 
\begin{align}
\frac{1+zx}{(1+zy_1)(1+zy_2)}&=
%z^{-1}(A(y_1)^{-1}+B(y_2)^{-1})-z^{-2}(A(y_1)^{-2}+B(y_2)^{-2})+z^{-3}(A(y_1)^{-3}+B(y_2)^{-3})-z^{-4}(A(y_1)^{-4}+B(y_2)^{-4})+...\\
%\frac{A}{1+zy_1}+\frac{B}{1+zy_2}&=
%\sum_{i\leq -1}z^{i}(-1)^{i+1}(A(y_1)^{i}+B(y_2)^{i})=
\sum_{l\leq -1}z^{l}c_{l}^{*} \, , \\
c_l^* &= -(-1)^{l}(A(y_1)^{l}+B(y_2)^{l}) \, , \ l \leq -1. \label{eq:ck*_closed_formula}
\end{align}
Such coefficients $c^*_l$ satisfy $c_{l}^{*}+c_{l-1}^{*}(y_1+y_2)+c_{l-2}^{*}y_1y_2=0$. %After shifting the index, $k=l-2$, the relation can be rewritten as
%\begin{align}
%& c_k^* y_1y_2 + c_{k+1}^* (y_1+y_2) + c_{k+2}^* = 0  \, , \ k \leq -3 \, ,  \\
%&c_{-1}^* = \frac{x}{y_1y_2} \, , \nonumber \\
%&c_{-2}^* = \frac{1}{y_1y_2}-\frac{x(y_2+y_1)}{y_1^{2}y_2^{2}} \, . \nonumber
%\end{align}
%Keep in mind that now
%\begin{equation}
%\frac{1+zx}{(1+zy_1)(1+zy_2)}=\sum_{k\leq -3}z^{k+2}c_{k+2}^{*} \, . 
%\end{equation}
The case not covered by Theorem~\ref{thm:REC-REL_ZERO} is the expansion in the region $\frac{1}{|y_1|}<|z|<\frac{1}{|y_2|}$, which is
\begin{align}
	\label{eq.annul}
	\frac{1+zx}{(1+zy_1)(1+zy_2)}
	%=A+B-z(A y_1+B y_2)+z^2(A(y_1)^2+B(y_2)^2)-z^3(A(y_1)^3+B(y_2)^3)+\dots \\    
	&
	%=\sum_{k=o}^{\infty}(z)^{k}(-1)^{k}(A(y_1)^k+B(y_2)^k)
	=\sum_{k=0}^{\infty}(-1)^{k} A(y_1)^k z^{k} + \sum_{l\leq -1}-(-1)^{l} B(y_2)^{l} z^{l}  \, .	
\end{align}
%
%The coefficient at $z$ is then 
%\begin{equation}
%	A y_1 = \frac{(x-y_1)y_1}{y_2-y_1}
%\end{equation}
\end{example}

Our goal will be to explore the expansions in the annular regions between the poles like~\eqref{eq.annul} --- not just in the above example, but in the general case.
In dimension $n|m$, the characteristic function $\Ber(E + zA)$ has generically $m$ different poles and we aim to understand the geometric meaning of the coefficients of these expansions.

One suggestive coincidence is that the number of expansion regions (including the neighborhoods of zero and infinity) is $m+1$, which matches the number of the Voronov-Zorich spaces (see the definition in the next section):
\begin{equation*}
	\Lambda^{k|0}(V) = \Lambda^k(V) \, , \quad \dots \, , \quad \Lambda^{k|m}(V)\, .  
\end{equation*}
These space  first appeared in the differential-geometric setup, which we  will  recall now. We will come back to the expansions in the algebraic setting in Section~\ref{sec.more}.

\section{Voronov-Zorich  $r|s$-forms}
\label{Sec:VZ}
An analogue of differential forms for supermanifolds that are graded by super-dimensions --- known as
$r|s$-forms --- were introduced and developed by Voronov-Zorich~\cite{tv:compl,tv:pdf,tv:cohom-1988} and Voronov~\cite{tv:git,tv:dual}. Here we recall them in the form suitable for our purposes\,\footnote{We do not touch $r|s$-forms with   $r<0$ introduced in~\cite{tv:dual}, though they are necessary in general picture.}.

On a supermanifold $M$ of dimension $n|m$, an  {$r|s$-form}  is a function $L=L(x,\dot x)$ of a point and $r$ even and $s$ odd tangent vectors, 
\[
L  \in C^\infty \left( \underbrace{TM \times_M \dots \times_M TM}_{r \ \text{times}} \times_M \underbrace{\Pi TM \times_M \dots \times_M \Pi TM}_{s \ \text{times}}  \right) \, ,
\]
satisfying conditions outlined below that  generalize the properties of the classical differential forms on ordinary manifolds. (A function $L$ is typically defined not everywhere in fiber arguments, but  has a singularity.) 
Here, $x=(x^A)$ represents local coordinates of a point on $M$, and $\dot x=(\dot 
x^A_F)$ collectively represents the components of $r$ even and $s$ odd tangent vectors at this point.
%The properties to satisfy are:$L(x, \dot x \cdot g)=L(x, \dot x)\cdot\Ber g$ 
% for all $g\in GL(r|s)$, and 
% \begin{equation*}
%	\frac{\partial^2 L}{\partial \dot x^A_F \partial \dot x^B_G} + 
%	(-1)^{\tilde{F} \tilde{G}+\tilde B (\tilde{F} + \tilde{G})} 
%	\frac{\partial^2 L}{\partial \dot x^A_G \partial \dot x^B_F} = 0\,,
%\end{equation*}
%where $\dot x_F^A$ are components of $x_F$, $A \in \{1, \dots, n, \hat 1, 
%\dots, \hat m \}$.

It is more convenient here to use an algebraic reformulation. 

If we denote $V=T_x M$ and $V^*=T^*_x M$, then $r|s$-forms (at a point $x$)
will be by definition the elements of $\Lambda^{r|s}(V^*)$. %(which is the  definition of the latter vector space). 
(Henceforth we drop mentioning of a point of $M$.) 
 We can think of $V$ as an arbitrary vector superspace of $\dim V=n|m$. Then \emph{$r|s$-forms} on  $V$ are defined as functions $L=L(\dot x)$, whose arguments are even and odd elements of $V$, satisfying the conditions:
\begin{equation}\label{eq.bercond} \tag{BER}
    L(g \cdot \dot x)=L( \dot x)\cdot\Ber g\,,
\end{equation}
for all $g\in GL(r|s)$, and 
\begin{equation}\label{eq.fund} \tag{PDE}
 \frac{\partial^2 L}{\partial \dot x^A_F \partial  \dot x^B_G} + 
(-1)^{\tilde{F} \tilde{G}+\tilde B (\tilde{F} + \tilde{G})} 
\frac{\partial^2 L}{\partial \dot x^A_G \partial  \dot x^B_F} = 0\,,
\end{equation}
for all indices, where $\dot x_F^A$ are components of $\dot x_F$, $A \in \{1, \dots, n, \hat 1, 
\dots, \hat m \}$.

By swapping the roles of $V$ and $V^*$, one   similarly defines the space  $\Lambda^{r|s}(V)$ of  \emph{$r|s$-vectors}  (or multivectors of degree $r|s$)  associated with a vector superspace $V$. See in~\cite{shemya:voronov2019:superplucker}. Properties of $\Lambda^{r|s}(V)$ and $\Lambda^{r|s}(V^*)$ are essentially the same (and differ only by a choice of $V$ or $V^*$ as a reference point). So, we will use one or another depending on convenience.

\begin{remark}
Conditions~\eqref{eq.bercond} and~\eqref{eq.fund} have the following differential-geometric origins. Condition~\eqref{eq.bercond} makes it possible to define an integral of an $r|s$-form on $M$ over an oriented $r|s$-path $\Gamma\,: U^{r|s}\to M$ so that it will be independent of parametrization. It extends to the supercase the notion of a ``$k$-density in an $n$-dimensional space'' (see e.g.~\cite{gelfand-gindikin-graev:intgeometry1980}) as was suggested in~\cite{Baranov_Shwarz}. 
Condition~\eqref{eq.fund} was found by Voronov-Zorich as necessary for defining an analogue of de Rham differential~\cite{tv:compl}, see explanations in~\cite{tv:git}. 
\end{remark}

Direct computations verify the following   lemmas.
\begin{lemma} 
Swapping indices $F$ and $G$ in equation~\eqref{eq.fund} yields an 
equivalent equation.
\end{lemma}
\begin{lemma} 
Swapping indices $A$ and $B$ in equation~\eqref{eq.fund} yields an 
equivalent equation.
\end{lemma}

One can also prove the following.

\begin{lemma}
	\label{lem:lin}
	 Every $r|s$-form is multilinear in its 
even arguments.
\end{lemma}
\begin{proof}
Considering \eqref{eq.fund} for $F=G$ and $F$ even, 
we get that $L$ is affine in even arguments.
Applying \eqref{eq.bercond} for $g\in GL(r|0)\subset GL(r|s)$, we in particular see   that $L$ must be homogeneous of degree $+1$ in each even argument. Hence it is  
linear in  even arguments.  
\end{proof}

Suppose $s=0$, i.e.   consider $r|0$-forms $L$ (where   the arguments are  even vectors). 

\begin{corollary}
   There is an isomorphism between the spaces  $\Lambda^{r|0}(V^*)$  and  $\Lambda^r(V^*)$ (the $r$th exterior power of $V^*$).   
\end{corollary}
\begin{proof}
We consider  $\Lambda^r(V^*)$   as the space of all antisymmetric $r$-linear functions on $V$. 
  From Lemma~\ref{lem:lin}, an $r|0$-form $L$ is a multilinear function of its  arguments (even vectors). By multilinearity, it can be extended   to  a function of $r$ vectors of arbitrary parity. Antisymmetry can be obtained either from~\eqref{eq.bercond} or from writing $L$ explicitly as a multilinear form and applying \eqref{eq.fund}. Conversely, if $L$ is an antisymmetric multilinear function, it will clearly satisfy~\eqref{eq.bercond}, and condition~\eqref{eq.fund} can be checked directly. 
\end{proof}

In particular, $r|0$-forms do not vanish for $r>n$ unless $m=0$. 

Lemma~\ref{lem:lin} implies that every $r|s$-form can be written as:
\begin{equation} \label{eq:L_gen}
L = \dot x_1^{A_1} \dots \dot x_r^{A_r} L_{A_r \dots A_1}(\dot x_{\hat 1}, \dots, \dot x_{\hat m}) \, ,
\end{equation}
 where $L_{A_r \dots A_1}$ is antisymmetric in its 
indices. 
%  Let us now consider all the remaining conditions that arise 
% from~\eqref{eq.fund}.
% \begin{enumerate}
%  \item Case: Both $x_F=x_k$ and $x_G=x_l$ are even. 
% Then: 
% \begin{equation}\label{eq.fund001} \tag{PDE00}
%  \frac{\partial^2 L}{\partial x^A_k \partial  x^B_l} = - 
% \frac{\partial^2 L}{\partial  x^A_l \partial  x^B_k} \,.
% \end{equation}
%  \item Case: Both $x_F=x_\mu$ and $x_G=x_\nu$ are odd. 
% Then: 
% \begin{equation}\label{eq.fund11} \tag{PDE11}
%  \frac{\partial^2 L}{\partial x^A_\mu \partial  x^B_\nu} = 
% \frac{\partial^2 L}{\partial  x^A_\nu \partial  x^B_\mu} \,.
% \end{equation}
%  \item Case: $x_F=x_k$ is even and $x_G=x_\mu$ is odd. 
% Then: 
% \begin{equation}\label{eq.fund10} \tag{PDE10}
%  \frac{\partial^2 L}{\partial  x^B_\mu \partial x^A_k } = 
% \frac{\partial^2 L}{\partial  x^A_\mu \partial  x^B_k} 
% (-1)^{(\tilde{A}+1)(\tilde{B}+1)}\,.
% \end{equation}
%  \item Case: Both $x_F$ and $x_G$ are odd. 
% Then we have a trivial identity.
% \end{enumerate}

\begin{lemma} \label{lem.homogen}
$L$ and each $L_{A_r \dots A_1}$ are homogeneous functions of 
	$x_{\hat 1}, \dots, x_{\hat m}$ of degree $-1$.
\end{lemma}
\begin{proof} For $L$, it immediately follows from~\eqref{eq.bercond}. Taking into account~\eqref{eq:L_gen}, we see that it holds for each coefficient
	$L_{A_r \dots A_1}$. 
\end{proof}

From   Lemma~\ref{lem.homogen}, it follows that $s$ odd vectors in the arguments of an $r|s$-form cannot be zero. Moreover, condition~\eqref{eq.bercond} implies that these vectors cannot be linearly dependent. This leads to the important conclusion that on an $n|m$-dimensional space, $s$ must satisfy
\begin{equation}
  0 \leq s\leq m\,.
\end{equation}
At the same time, we saw that there are no restrictions on $r$. 

\begin{lemma} Equality~\eqref{eq.bercond} implies
\begin{equation} \label{eq:fromBer} \tag{BER'}
  \dot x_G^A \frac{\partial L}{\partial \dot x_F^A}- (-1)^{\tilde{F}} 
\delta^F_G L =0 \, ,
\end{equation}
for all $F,G,A,B$. Here $\delta^F_G=1$ for $F=G$ and is zero otherwise.
\end{lemma}
\begin{proof}
 Plug in $g$ which is close to the identity.
\end{proof}

\begin{remark}
	The expression 
\begin{equation}
T^G_F= \dot x_G^A \frac{\partial L}{\partial \dot x_F^A}- 
(-1)^{\tilde{F}} 
\delta^F_G L 
\end{equation}
is well-known in physics as the \emph{energy-momentum tensor} or 
\emph{stress tensor}. Equation~\eqref{eq:fromBer} gives the infinitesimal form of the homogeneity properties of $L$ and the invariance of $L$ under elementary transformations (adding a multiple of one argument to another argument). Thus it is  equivalent to~\eqref{eq.bercond}  for $g$ in the component of identity in $GL(r|s)$.   
\end{remark}

Suppose $s=m$. By using~\eqref{eq.bercond}, one can show that $r|m$-forms vanish for $r>n$ (in contrast with $r|0$-forms). Moreover, the following descriptions exist for the ``top'' case of $n|m$-forms and for $r|m$-forms with $r<n$.

\begin{theorem}[\cite{tv:git}]  \label{thm:Ber}
	The ``top'' Voronov-Zorich forms are Berezinians (up to a multiplication by a constant).
	Namely, Voronov-Zorich $n|m$-forms on an $n|m$-dimensional space have the form:
	\begin{equation}
		L(\dot X)= C \Ber(\dot X) \, ,
	\end{equation}	
	where we use matrix re-formulation and 
	where $\dot X$ is the matrix formed by the arguments $\dot x$ of $L$ as rows. Here, $C$ is a constant.
\end{theorem}
\begin{proof}
	Consider~\eqref{eq.bercond} in the matrix form:
	$L(g\dot X) = L(\dot X) \Ber g$.
	Setting $g=\dot X^{-1}$, we get $L(\dot X)= L(1) \Ber(\dot X)$. One can   show (though it is not easy) that the function $\Ber(\dot X)$ satisfies also condition~\eqref{eq.fund}.	
\end{proof}

\begin{theorem}[\cite{tv:git}] \label{thm:rmforms}
	Every Voronov-Zorich $r|m$-form on an $n|m$-dimensional space  has the following form:
	\begin{equation}
		L(\dot x_1,\dots, \dot x_r, \dot x_{\hat 1}, \dots ,\dot x_{\hat m}) =
		\Ber 
		\begin{pmatrix}
			K & \mathcal{M} \\
			K' & \mathcal{M}' \\
			\mathcal{L} & N 
		\end{pmatrix}		\, ,
	\end{equation}
	where the blocks $K$ and $\mathcal{M}$ correspond to the even arguments $\dot x_1,\dots, \dot x_r$ written as rows, while the blocks $\mathcal{L}$ and $N$ correspond to the odd arguments $\dot x_{\hat 1}, \dots ,\dot x_{\hat m}$ also written as rows. The blocks $K'$ and $\mathcal{M}'$ provide additional $n-r$ even rows in the matrix so to make it square. They serve as coefficients and may be of any parity.
\end{theorem}

(Theorem~\ref{thm:rmforms} can be reformulated as establishing an isomorphism between $r|m$-forms and ``integral forms'',  which we recall in Section~\ref{sec:BerLei}.)  

While quite a lot is known about $r|s$-forms in general,  particularly their differential-geometric properties as part of integration theory on supermanifolds, and there are various  examples  including estimates from below for their cohomology, see~\cite{tv:cohom-1988} and \cite{tv:git}, there is no explicit description of $r|s$-forms apart from the two special cases $s=0$ and $s=m$. This remains a challenging problem. 

In the general case, the spaces $\Lambda^{r|s}(V^*)$ should be infinite-dimensional --- from known examples containing functional parameters~\cite{tv:git}. (See also our   result in the next section concerning  $1|1$-forms.)   
So, one may naturally wish to make these spaces ``smaller'' by suitable choices. Another aspect to be taken into account is   ambiguity in the domain of     definition of an $r|s$-form. Since its $s$ odd arguments  must be linearly independent, should the  form be defined everywhere subject only to this restriction? Original publications~\cite{tv:compl,tv:pdf,tv:cohom-1988} and \cite{tv:git} seem to suggest just that. However, Belopolsky~\cite{belopolsky1997new} attracted attention specifically  to $r|s$-forms with singularities (they had been considered before~\cite{hov:gayduk}). (Similar examples of   $r|s$-vectors appeared in~\cite{shemya:voronov2019:superplucker}.)

The strategy that we want to follow is to identify a class of  $r|s$-forms possibly with singularities that, on one hand, is ``small'' enough so to admit an explicit description, and, on another hand, is ``large'' enough for all purposes. In particular it should be closed  under the exterior and interior product by (co)vectors (i.e. the action of the Clifford algebra of $V\oplus V^*$) and in the differential-geometric context closed under differential. It should also be functorial. We have not discussed functorial properties of  $r|s$-forms yet; Voronov and Zorich identified the class of \emph{proper morphisms} of supermanifolds for this purpose. In algebraic setting, this means linear maps $V\to W$  monomorphic in the odd directions;  a reasonable class of $r|s$-forms must admit pullbacks with respect to all such morphisms.

\section{$1|1$-forms  as differential forms on the projective superspace}
\label{Sec:11}

Consider the case of $1|1$-forms at a point $x\in M$, or,   in the algebraic language, the space  $\Lambda^{1|1}(V^*)$ for a vector superspace $V$. 
Such forms are functions  of one even 
vector $x_1$ and one odd vector $x_{\hat 1}$.
Here we drop the dot and switch our notations from $\dot x$ to $x$. (Since $x$, previously denoting    a point of a supermanifold, is no longer in use.)

\subsection{Detalization of properties and some examples}
From the previous analysis we know that: 
%(which also depend on $x$, though we omit it for brevity): 
\begin{equation} \label{eq:L11}
L = 
x_1^A L_A \, ,
\end{equation}
where $L_A=L_A(x_{\hat 1})$, and $A \in \{1,\dots,n,\hat{1}, \dots ,\hat{m} \}$.

Let us now consider conditions that arise 
from~\eqref{eq.fund}.
\begin{enumerate}
 \item Cases: $\dot x_F= \dot x_G=x_1$ or $\dot x_F= \dot x_G=x_{\hat 1}$ results in a
trivial identity. 
 \item Case: $\dot x_F=x_1$ is even and $\dot x_G=x_{\hat 1}$ is odd. 
Then: 
\begin{equation}\label{eq:PDE'} \tag{PDE'}
 \frac{\partial L_A}{\partial  x^B_{\hat 1}} = 
\frac{\partial L_B}{\partial  x^A_{\hat 1}} 
(-1)^{(\tilde{A}+1)(\tilde{B}+1)}\,.
\end{equation}
\end{enumerate}

Consider condition~\eqref{eq:fromBer}. Setting $\dot x_F = \dot 
x_G = 
x_1$, we obtain $L = 
x_1^A L_A$. Setting $\dot x_F = x_1$ and $\dot 
x_G = 
x_{\hat 1}$, we obtain
\begin{equation} \label{eq:BER1'} \tag{BER1'}
x^{A}_{\hat 1} L_{A}=0 \, . 
\end{equation}

Setting $\dot x_F = 
x_{\hat 
1}$ and 
$\dot 
x_G =x_1$ yields an identity that is trivial modulo the conditions 
already 
obtained.
% 
%  $x_1^A \frac{\partial L}{\partial x_{\hat 1}^A}=0$,
%  which implies 
%  $x_1^1 \frac{\partial L}{\partial x_{\hat 1}^1}+
%  x_1^{\hat 1} \frac{\partial L}{\partial x_{\hat 1}^{\hat 1}}=0$, 
%  which can be re-written as
%  $x_1^1 x_1^1 \frac{\partial L_1}{\partial x_{\hat 1}^1}+
%  x_1^{\hat 1} x_1^{\hat 1} \frac{\partial L_{\hat 1}}{\partial 
%  x_{\hat 
%  1}^{\hat 1}}=0$, implying 
% $\frac{\partial L_1}{\partial x_{\hat 
%  1}^1}=0$. 
Finally, writing 
condition~\eqref{eq:fromBer} with $\dot x_F = x_{\hat 
1}$ and 
$\dot 
x_G =x_{\hat 1}$ yields
\begin{equation} \label{eq:BER2'} \tag{BER2'}
x_{\hat 1}^A \frac{\partial L}{\partial 
x_{\hat 1}^A} =-L \, . 
\end{equation}

%The ultimate goal now is to explicitly describe the 
%solutions to the 
%system of 
%equations~\eqref{eq:L11},~\eqref{eq:PDE'},\eqref{eq:BER1'},
%~\eqref{eq:BER2'}.

\begin{example}[An illustration for Theorem~\ref{thm:Ber}] \label{ex11in11}
	Consider $1|1$-forms in a 
	$1|1$-dimensional space. If the 
	components are denoted as 
	$x_1 = (x^1_1 | x^{\hat 1}_1)$ and 
	$x_{\hat 1} = (x_{\hat 1}^1 | x_{\hat 1}^{\hat 1})$, we have a 
	linear system (\eqref{eq:L11} and ~\eqref{eq:BER1'}):
	\begin{equation}
		\begin{bmatrix}
			\begin{array}{c|c}
				x_1^1 & x_1^{\hat 1} \\
				\hline
				x_{\hat 1}^1 & x_{\hat 
					1}^{\hat 1} \end{array}
		\end{bmatrix}
		\begin{bmatrix} L_{1} \\ \hline L_{\hat 1}
		\end{bmatrix}
		=
		\begin{bmatrix} L \\ \hline 0 \end{bmatrix} \, .
	\end{equation}
	Using the super analog of Cramer's rule~\cite{Kh_V_recurrent} (see 
	also~\cite{shemya:voronov2019:superplucker} for more details), we obtain:
	\begin{align}
		L_1 & = \frac{\Ber 
			\begin{bmatrix}
				\begin{array}{c|c}
					L & x_1^{\hat 1} \\
					\hline
					0 & x_{\hat 
						1}^{\hat 1} \end{array}
		\end{bmatrix}}{\Ber X} 
		=\frac{L/x_{\hat 
				1}^{\hat 1}}{\Ber X}
		=\frac{L}{x^1_1-x^{\hat 1}_1(x_{\hat 1}^{\hat 1})^{-1} x_{\hat 1}^1}
		\, , \\
		%%%%%%%%%%%%%%%%%%%%%%%%%%%%%%%%%%%%%%%%
		L_{\hat 1} & = \frac{\Ber^* 
			\begin{bmatrix}
				\begin{array}{c|c}
					x_1^1 & L \\
					\hline
					x_{\hat 1}^1 & 0
				\end{array}
		\end{bmatrix}}{\Ber^*  X} =
		\frac{\Ber
			\begin{bmatrix}
				\begin{array}{c|c}
					0 & x_{\hat 1}^1 \\
					\hline
					L & x_1^1
				\end{array}
		\end{bmatrix}}{\Ber 
			\begin{bmatrix}
				\begin{array}{c|c}
					x_{\hat 
						1}^{\hat 1} & x_{\hat 1}^1 \\
					\hline
					x_1^{\hat 1} & x_1^1  \end{array}
			\end{bmatrix}
		} 
		=
		\frac{-x_{\hat 1}^1 (x_1^1)^{-1} L}{x_{\hat 
				1}^{\hat 1} - x_{\hat 1}^1 (x_1^1)^{-1}x_1^{\hat 1}} =
		\frac{-x_{\hat 1}^1  L}{x_{\hat 
				1}^{\hat 1} x_1^1 - x_{\hat 1}^1 x_1^{\hat 1}}
		=
		\frac{-x_{\hat 1}^1  L}{x_{\hat 
				1}^{\hat 1} x_1^1 }
		\, ,
	\end{align}
	with $X$ being the matrix formed by the arguments $x_1$ and $x_{\hat 
		1}$ as rows (i.e., the coefficient matrix of the above linear 
	system). 
	
	Now, altogether  conditions \eqref{eq:L11}, \eqref{eq:PDE'}, 
	\eqref{eq:BER1'},
	\eqref{eq:BER2'} are equivalent to the following:
	\begin{align} \label{eq:Lgenform}
		& L = 
		x^1_1 L_1 + x^{\hat 1}_1 L_{\hat 1}
		\, , \\ \label{eq:L1Lhat1viaL}
		& 
				L_1 =\frac{1}{x^1_1-x^{\hat 1}_1(x_{\hat 1}^{\hat 1})^{-1} 
			x_{\hat 1}^1} L
		\, , \
		L_{\hat 1} =-
		\frac{x_{\hat 1}^1}{x_{\hat 
				1}^{\hat 1} x_1^1 } L
		\, , \\
		& \frac{\partial L_1}{\partial x_{\hat 1}^1} = 0 \, ,\  
		\frac{\partial L_{1}}{ \partial x_{\hat 1}^{\hat 1} } =  
		\frac{\partial L_{\hat 1}}{ \partial x_{\hat 1}^1 } \, , \\
		& x_{\hat 1}^1 \frac{\partial L}{\partial 
			x_{\hat 1}^1} + x_{\hat 1}^{\hat 1} \frac{\partial L}{\partial 
			x_{\hat 1}^{\hat 1}} =-L \, . \label{eq:der}
	\end{align}
	From~\eqref{eq:L1Lhat1viaL}, we obtain $L=(x^1_1-x^{\hat 1}_1(x_{\hat 
		1}^{\hat 1})^{-1} 
	x_{\hat 1}^1) L_1$, where $L_1$ does not depend on $x^1_{\hat 1}$. 
	Substituting this into~\eqref{eq:der}, we obtain:
	\begin{equation}
		x_{\hat 1}^{\hat 1} \frac{\partial L_1}{\partial 
			x_{\hat 1}^{\hat 1}} = - L_1 \, ,
	\end{equation}
	and so
	\begin{equation}
		L_1 = C (x_{\hat 1}^{\hat 1})^{-1} \, ,
	\end{equation}
	where $C$ is a constant (does no depend on arguments $x_1$ and 
	$x_{\hat 1}$).
	Hence, 
	\begin{equation} \label{11in11}
		L(x_1|x_{\hat 1}) = C \Ber X = C \frac{x^1_1-x_1^{\hat 1} (x_{\hat 
				1}^{\hat 1})^{-1} x_{\hat 1}^1}{x_{\hat 1}^{\hat 1}}\, .
	\end{equation}
	One can verify that~\eqref{11in11} indeed satisfies 
	conditions~\eqref{eq:Lgenform}-\eqref{eq:der}.
\end{example}

\begin{example}[An illustration for Theorem~\ref{thm:rmforms}]  \label{ex:11inn1} Consider the case of $1|1$-forms on the 
	$n|1$-dimensional space: 
	\begin{equation}
		L =L (x_1, x_{\hat 1}) =  x^1_1 L_1(x_{\hat 1} ) + \dots + 
		x^n_1 L_1(x_{\hat 1} ) + x^{\hat 1}_1 L_{\hat 1}(x_{\hat 1} )\, . 
	\end{equation}
	
	In this case, conditions~\eqref{eq.bercond} and~\eqref{eq.fund} become:
	\begin{align}
		\label{eq:LL}& L = 
		x^i_1 L_i + x^{\hat 1}_1 L_{\hat 1}
		\, , \\
		& 
		\frac{\partial L_j}{\partial x_{\hat 1}^i} = - \frac{\partial 
			L_i}{\partial x_{\hat 1}^j} \, , \
		\frac{\partial 
			L_i}{\partial x_{\hat 1}^i} = 0 \, , \
		\frac{\partial L_{i}}{ \partial x_{\hat 1}^{\hat 1} } =  
		\frac{\partial L_{\hat 1}}{ \partial x_{\hat 1}^i } \, , \\
		& \label{eq:iszero}
		x^i_{\hat 1} L_i + x^{\hat 1}_{\hat 1} L_{\hat 1}=0 \, , \\
		& x_{\hat 1}^{i} 
		\frac{\partial L}{\partial 
			x_{\hat 1}^i}+ x_{\hat 1}^{\hat 1} 
		\frac{\partial L}{\partial 
			x_{\hat 1}^{\hat 1}}=-L \, , \label{eq:minusL}
	\end{align}
	where $i,j \in \{1, \dots, n \}$.
	%  and all $a \in \{1, \dots, n, \hat 1 \}$.
		Expressing $L_{\hat 1}$ from~\eqref{eq:iszero}, and substituting it 
	into~\eqref{eq:LL}, we obtain:
	\begin{equation} \label{eq:LLi}
		L = (x_1^i - x_1^{\hat 1} (x^{\hat 1}_{\hat 1})^{-1} x_{\hat 1}^i) 
		L_i \, .
	\end{equation}
	Then, since $L_i$ does not depend on $x_{\hat 1}^i$,~\eqref{eq:LLi} implies:
	\begin{equation} \label{eq:LLi1}
		x_{\hat 1}^{i} 
		\frac{\partial L}{\partial 
			x_{\hat 1}^i} =  x_{\hat 1}^{i} x_1^{\hat 1} (x^{\hat 1}_{\hat 
			1})^{-1} L_i \, .
	\end{equation}
	Also,
	\begin{equation} \label{eq:LLi2}
		x_{\hat 1}^{\hat 1} 
		\frac{\partial L}{\partial 
			x_{\hat 1}^{\hat 1}} =  x_{\hat 1}^{\hat 1} \left(
		x_1^{\hat 1} (x^{\hat 1}_{\hat 1})^{-2} x_{\hat 1}^i L_i + 
		(x_1^i - x_1^{\hat 1} (x^{\hat 1}_{\hat 1})^{-1} x_{\hat 1}^i) 
		\frac{\partial L_i}{\partial x_{\hat 1}^{\hat 1}}  
		\right) 
	\end{equation}
	Substituting~\eqref{eq:LLi}, \eqref{eq:LLi1}, and \eqref{eq:LLi2} 
	into~\eqref{eq:minusL}, two terms cancel out, and we obtain: 
	% \begin{equation}
		%   x_{\hat 1}^{\hat 1}  
		% (x_1^i - x_1^{\hat 1} (x^{\hat 1}_{\hat 1})^{-1} x_{\hat 1}^i) 
		% \frac{\partial L_i}{\partial x_{\hat 1}^{\hat 1}}  
		%   = - (x_1^i - x_1^{\hat 1} (x^{\hat 1}_{\hat 1})^{-1} x_{\hat 1}^i) 
		% L_i  \, .
		% \end{equation}
	% This implies: 
	\begin{equation} \label{firstfun}
		(x_1^i - x_1^{\hat 1} (x^{\hat 1}_{\hat 1})^{-1} x_{\hat 1}^i) \left(
		x_{\hat 1}^{\hat 1}  
		\frac{\partial L_i}{\partial x_{\hat 1}^{\hat 1}}  
		+
		L_i \right) =0 \, .
	\end{equation}
	
	Since $L_i$ is independent of $x_1$, for every $i \in \{1, \dots , 
	n\}$, we can set $x_1$ in~\eqref{firstfun} step by step to be 
	equal basis vectors: $x_1^j = \delta_i^j$, and $x_1^{\hat 1}=0$. This 
	leads to the second factors in the sum~\eqref{firstfun} being 
	zeros for all values of $i$. This implies that the general solution 
	to~\eqref{firstfun} is:
	\begin{equation}
		L_i = \frac{C_i}{x^{\hat 1}_{\hat 1}} 
	\end{equation}
	for all $i \in \{1, \dots, n \}$. Substituting this 
	into~\eqref{eq:LLi}, we obtain:
	\begin{equation}
		L = C_i \Ber X^{i \hat 1} \, , 
	\end{equation}
	where $X^{i \hat 1}$ is the submatrix formed by the 
	$i$th and $\hat 1$th columns of
	\begin{equation}
		X=
		\begin{bmatrix}
			\begin{array}{c|c}
				x_1^1 \dots x_1^n & x_1^{\hat 1} \\
				\hline
				x_{\hat 
					1}^1 \dots x_{\hat 
					1}^n & x_{\hat 
					1}^{\hat 1} \end{array}
		\end{bmatrix} \, .
	\end{equation}
\end{example}

\subsection{Relation with the projective superspace}

We propose the following interpretation of $1|1$-forms\,\footnote{The possibility of a relation between $1|1$-forms and forms on a projective space was suggested by A.~Odesskii in a private communication.}. Consider 
\begin{equation}
(x_{\hat 1}^{\hat 1}: \ldots : x_{\hat  1}^{\hat m}: x_{\hat 1}^1: \ldots : x_{\hat 1}^n) 
\end{equation}
(note the order) as a point of the projective superspace $\mathbb{P}^{m-1|n}$, so that the components of the odd vector $x_{\hat 1}$ are regarded as homogeneous coordinates. We write them swapping the order so  to obtain the standard format in which even variables go first and odd variables go second. We associate   the components of the even vector $x_1$ with the differentials of the components of  $x_{\hat 1}$\,:
\begin{equation}
x_1\sim (d x_{\hat 1}^1,\, \ldots\,,  d x_{\hat 1}^n, dx_{\hat 
1}^{\hat 1},\, \ldots\,, dx_{\hat 
1}^{\hat m}) \,. 
\end{equation}
(We view   $d$ as an odd operation.)
In this way, we consider the following differential $1$-form on an $m|n$-dimensional superspace $\Pi V$ as associated 
with a $1|1$-form~\eqref{eq:L11}  on $V$\,:
\begin{equation}
\omega_L:=dx_{\hat 1}^a 
L_a \, . 
\end{equation}
Coordinates on $\Pi V$ are homogeneous coordinates from the viewpoint of $\mathbb{P}^{m-1|n}$. Our plan is to relate $\omega_L$ with a form on the projective superspace $\mathbb{P}^{m-1|n}$. This goes as follows.

Consider the chart for $\mathbb{P}^{m-1|n}$ where the 
homogeneous coordinate
$x_{\hat 1}^\alpha$ is invertible, here $\alpha$ can be any from $\{\hat{1}, \dots, \hat{m}\}$. Denote the corresponding inhomogeneous 
coordinates by $u_{\hat 1}^a$, 
\begin{equation}
x_{\hat 1}^{\alpha} u_{\hat 1}^a = x_{\hat 1}^a \, 
\end{equation}
for each $a \in \{1, \dots, n, \hat 1, \dots, \hat m 
\}$. Here $u_{\hat 1}^\alpha=1$ is a fake coordinate.

By $\overline{L_a}$ we denote
$L_a$ with   $u_{\hat 1}^a$ substituted for $x_{\hat 1}^a$ %inhomogeneous coordinates: 
%\[
%u_{\hat{1}}^i = \frac{x_{\hat{1}}^i}{x_{\hat{1}}^\alpha}, \quad
%u_{\hat{1}}^{\hat{\jmath}} = \frac{x_{\hat{1}}^{\hat{\jmath}}}{x_{\hat{1}}^\alpha},
%%\]
\begin{equation} 
	\overline{L_a} = L_a(u_{\hat{1}}^1, \dots, u_{\hat{1}}^n \mid u_{\hat{1}}^{\hat{1}}, \dots, u_{\hat{1}}^{\hat{m}}) \, , \ u^\alpha_{\hat 1}=1 \, . 
\end{equation}

\begin{lemma} 
	For $a \in \{1, \dots, n, \hat 1, \dots, \hat m 
	\}$, we have:
	\begin{equation} \label{eq:BER1''11innm}
		\tag{BER1''}
		u_{\hat 1}^a \overline{L_a} = 0\, .  
	\end{equation}
\end{lemma}
\begin{proof}
		Condition~\eqref{eq:BER1'} implies $0 = x_{\hat 1}^a L_a = 
	x_{\hat 1}^{\alpha} u_{\hat 1}^a L_a = u_{\hat 
		1}^a \overline{L_a}$, and so we have~\eqref{eq:BER1''11innm}.
\end{proof}

\begin{lemma}
	\begin{equation}\label{eq:PDE''} \tag{PDE''}
		\frac{\partial \overline{L_a}}{\partial  u^b_{\hat 1}} = 
		\frac{\partial \overline{L_b}}{\partial  u^a_{\hat 1}} 
		(-1)^{(\tilde{a}+1)(\tilde{b}+1)}\, ,
	\end{equation}	
	for each $a,b \in \{1, \dots, n, \hat 1, \dots, \hat m 
	\}$, $a \neq \alpha$, $b \neq \alpha$. 
\end{lemma}
\begin{proof}
	Since $L_a$ is homogeneous of degree $-1$, we have:
	\begin{multline*}
			L_a(
	x_{\hat{1}}^1, \dots, x_{\hat{1}}^n \mid x_{\hat{1}}^{\hat1}, \dots, x_{\hat{1}}^{\hat{m}})
	= L_a\left(  x_{\hat{1}}^\alpha \frac{x_{\hat{1}}^1}{x_{\hat{1}}^\alpha}, \dots, x_{\hat{1}}^\alpha \frac{x_{\hat{1}}^n}{x_{\hat{1}}^\alpha} \,\Big|\, x_{\hat{1}}^\alpha \frac{x_{\hat{1}}^{\hat{1}}}{x_{\hat{1}}^\alpha}, \dots, x_{\hat{1}}^\alpha \frac{x_{\hat{1}}^{\hat{m}}}{x_{\hat{1}}^\alpha} \right) \\
	= \left(x_{\hat{1}}^\alpha\right)^{-1} \cdot \overline{L_a} (
	u_{\hat{1}}^1, \dots, u_{\hat{1}}^n \mid u_{\hat{1}}^{\hat1}, \dots, u_{\hat{1}}^{\hat{m}}) \, .
	\end{multline*}
	That is:
\begin{equation} \label{eq:overLandL}
   x_{\hat 1}^\alpha L_a(x) =   \overline{L_a}(u) \, .
\end{equation}
%where $L_a=L_a(
%x_{\hat 1}^1,\dots,x_{\hat 1}^n|x_{\hat 1}^{\hat 1}, \dots, 
%x_{\hat 1}^{\hat m})$, and 
%$\overline{L_a}=L_a(
%u_{\hat 1}^1,\dots,u_{\hat 1}^n|u_{\hat 1}^{\hat 1}, \dots, 
%u_{\hat 1}^{\hat m})$.

Differentiating equality~\eqref{eq:overLandL}, we obtain:
\begin{equation} \label{eq:connectders}
 x_{\hat 1}^\alpha  \frac{\partial L_a}{\partial x_{\hat 1}^b} =  
  \frac{\partial \overline{L_a}}{\partial x_{\hat 1}^b} =
\frac{1}{x_{\hat 1}^\alpha}
\frac{\partial \overline{L_a}}{\partial u_{\hat 1}^b} 
\end{equation}
for each $a,b \in \{1, \dots, n, \hat 1, \dots, \hat m 
\}$, $a \neq \alpha$, $b \neq \alpha$. This implies that~\eqref{eq:PDE'} for $L$ is equivalent to~\eqref{eq:PDE''} for $\overline{L}$.
\end{proof}

As we know, not every differential form  written in homogeneous coordinates %$x^a_{\hat 1}$ 
descends to the projective space. 
%A generic form is not projectively well-defined and, upon descent, may depend on the choice of representative within an equivalence class, rather than on the point in projective space itself.
However, we will see that in our case, the conditions satisfied by $L$ ensure that for the $1$-form $\omega_L$  associated with $1|1$-forms $L$ this descent is possible. 
\begin{lemma}
 Every differential $1$-form
  $\omega_L$ associated with a $1|1$-form on an $n|m$-dimensional 
space descends
  to the projective space $\mathbb{P}^{m-1|n}$.  
\end{lemma}
\begin{proof} Consider the chart defined by selecting
homogeneous coordinates where
$x_{\hat 1}^\alpha=1$.
% First, we note that the homogeneity of each $L_a$ implies
% \begin{equation}
%   x_{\hat 1}^\alpha L_a =  \overline{L_a} \, ,
% \end{equation}
% where $L_a=L_a(
% x_{\hat 1}^1,\dots,x_{\hat 1}^n|x_{\hat 1}^{\hat 1}, \dots, 
% x_{\hat 1}^{\hat m})$, and 
% $\overline{L_a}=L_a(
% u_{\hat 1}^1,\dots,u_{\hat 1}^n|u_{\hat 1}^{\hat 1}, \dots, 
% u_{\hat 1}^{\hat m})$, $a \in \{1, \dots, n, \hat 1, \dots, \hat m 
% \}$, and where 
% \begin{equation}
% x_{\hat 1}^{\alpha} u_{\hat 1}^a = x_{\hat 1}^a \, .
% \end{equation}
We show that we can re-write $\omega_L=dx_{\hat 1}^a L_a$ in terms of $u_{\hat 1}^a$ and $d u_{\hat 1}^a$ only, which means that $\omega_L$ indeed represents a differential form on the projective space:
\begin{align*}
\omega_L = dx_{\hat 1}^a L_a(x) &= d(x_{\hat 1}^{\alpha} u_{\hat 1}^a) 
L_a(x) \\
&= 
dx_{\hat 1}^{\alpha} u_{\hat 1}^a 
L_a(x) + x_{\hat 1}^{\alpha} d u_{\hat 1}^a L_a(x)\\
&= 
dx_{\hat 1}^{\alpha} \frac{u_{\hat 1}^a}{x_{\hat 1}^\alpha} 
\overline{L_a}(u) +  d u_{\hat 1}^a \overline{L_a}(u)  \quad (\text{we used}~\eqref{eq:overLandL})\\
&= d u_{\hat 1}^a \overline{L_a}(u) 
\quad (\text{we used}~\eqref{eq:BER1''11innm})
\, .
\end{align*}
\end{proof}

Thus, we have the projective descent of $\omega_L$:
\begin{equation}
\overline{\omega_L}	= d u_{\hat 1}^a \overline{L_a}(u) \, .
\end{equation}

\begin{lemma} 	\label{lem:closed_is_equiv_to_PDE}	
	\begin{enumerate}
		\item $\omega_L$ is closed if and only if 
		~\eqref{eq:PDE'} holds;
		\item $\overline{\omega_L}$ is closed if and only if 
		~\eqref{eq:PDE''} holds.
	\end{enumerate}
\end{lemma}
\begin{proof} 
	This is just the explicit form of the condition of closedness of a $1$-form on a supermanifold. For concreteness, let us see the computation for $\omega_L$:
	\begin{align*}
		d \omega_L &= (-1)^{\tilde{a}} dx_{\hat 1}^a 
		d L_a \\
		&= (-1)^{\tilde{a}} dx_{\hat 1}^a  
		dx_{\hat 1}^b \frac{\partial L_a}{\partial x_{\hat 1}^b} \\
		& 
		= \frac{1}{2} \left((-1)^{\tilde{a}} dx_{\hat 1}^a  
		dx_{\hat 1}^b \frac{\partial L_a}{\partial x_{\hat 1}^b} + 
		(-1)^{\tilde{a}} dx_{\hat 1}^a  
		dx_{\hat 1}^b \frac{\partial L_a}{\partial x_{\hat 1}^b} \right) \quad (\text{we duplicated the term})\\
		&=
		\frac{1}{2} \left( (-1)^{\tilde{a}}  dx_{\hat 1}^a  
		dx_{\hat 1}^b \frac{\partial L_a}{\partial x_{\hat 1}^b} + 
		(-1)^{\tilde{b}}  dx_{\hat 
			1}^b  
		dx_{\hat 1}^a \frac{\partial L_b}{\partial x_{\hat 1}^a} \right)  \\
		&=
		\frac{1}{2}  dx_{\hat 1}^a  
		dx_{\hat 1}^b \left( (-1)^{\tilde{a}}  \frac{\partial L_a}{\partial 
			x_{\hat 1}^b} + 
		(-1)^{\tilde{b}} (-1)^{\tilde{a}\tilde{b}}    
		\frac{\partial L_b}{\partial x_{\hat 1}^a} \right)
	\end{align*}
	where we essentially symmetrized the coefficients.
	 We see
	that the form is closed if and only if condition~\eqref{eq:PDE'} 
	holds. 
%		Taking into account~\eqref{eq:connectders}, we also prove the statement of the lemma.
\end{proof}

Combining everything together, we arrive at the following theorem.
\begin{theorem} \label{thm.11}
	On an $n|m$-dimensional supermanifold $M$,  the vector space of Voronov-Zorich $1|1$-forms at a point $p$ is isomorphic to the  vector space of closed differential $1$-forms on the projective superspace $\mathbb{P}^{m-1|n}=P(\Pi V)$, where $V=T_pM$.
\end{theorem}

We have used here the differential-geometric language, but obviously it can be reformulated without mentioning $M$. 

\begin{theorem} \label{thm.11alg}
	For an $n|m$-dimensional super vector space $V$,  the   space $\Lambda^{1|1}(V^*)$ is isomorphic to the  vector space of closed differential $1$-forms on  $\mathbb{P}^{m-1|n}=P(\Pi V)$. 
\end{theorem}

It would be very interesting to see if there is a connection between this result and the relation between the theory of $r|s$-forms and Gelfand-Gindikin-Graev integral geometry~\cite{gelfand-gindikin-graev:intgeometry1980} discovered in~\cite{tv:cohom-1988}. 

\begin{example} Consider again case of $1|1$-forms in a 
$1|1$-dimensional space: 
\begin{equation*}
 \omega_L = dx_{\hat 1}^1 L_1 + dx_{\hat 1}^{\hat 1} L_{\hat 1} \, ,
\end{equation*}
where $L_1=L_1(x_{\hat 1}^1|x_{\hat 1}^{\hat 1})$ and 
$L_{\hat 1}=L_{\hat 1}(x_{\hat 1}^1 \Big|x_{\hat 1}^{\hat 1})$. From the viewpoint of the theorem, we have a degenerate case, where the projective space $\mathbb{P}^{0|1}$ coincides with a $0|1$-dimensional affine superspace, and has a single chart with one odd 
coordinate $
 u^1_{\hat 1} = x_{\hat 1}^1 (x_{\hat 1}^{\hat 1})^{-1}$.
We already know that (see~\eqref{11in11} in Example~\ref{ex11in11})
%
%Condition~\eqref{eq:PDE'} for $a=b=1$ implies 
%\begin{equation*}
% \frac{\partial L_1}{\partial x^1_{\hat 1}}= 0\, .
%\end{equation*}
%Using homogeneity, we conclude that 
%\begin{equation*}
% L_1 = \frac{C_1}{x^{\hat 1}_{\hat 1}} \, .
%\end{equation*}
%From~\eqref{eq:BER1'}, we have
%\begin{equation*}
% L_{\hat 1} = -\frac{x^1_{\hat 1}}{x^{\hat 1}_{\hat 
%1}} L_1 =  \frac{-x^1_{\hat 1} C_1}{(x^{\hat 1}_{\hat 
%1})^2} \, .
%\end{equation*}
%Hence, 
\begin{equation*}
 \omega_L   = dx_{\hat 1}^1 L_1 + dx_{\hat 1}^{\hat 1} L_{\hat 1}= dx_{\hat 1}^1 \frac{C_1}{x^{\hat 1}_{\hat 1}} - dx_{\hat 
 1}^{\hat 1} \frac{x^1_{\hat 1} C_1}{(x^{\hat 1}_{\hat 
 1})^2}   
%= \frac{dx_{\hat 1}^1 -dx_{\hat 
%1}^{\hat 1}(x^{\hat 1}_{\hat 
%1})^{-1} x^1_{\hat 1} }{x^{\hat 1}_{\hat 
%1}} C_1 
= \Ber \begin{pmatrix}
\begin{array}{c|c}
dx_{\hat 1}^1 & dx_{\hat 
1}^{\hat 1}  \\
\hline
x_{\hat 1}^1 & x^{\hat 1}_{\hat 
1}
\end{array}
\end{pmatrix} C_1\, .
\end{equation*}
Now,  we have:
\begin{equation}
	\overline{L_1}=L_1(u^1_{\hat 1}, 1)= C_1 \, .
\end{equation}
This implies:
\begin{equation}
	\overline{\omega_L} = 	
	d u^1_{\hat 1} C_1 \, ,
\end{equation}
where we use the fact that for the other differential we have: 
$du_{\hat 
	1}^{\hat 1}=0$.
\end{example}

\section{Bernstein-Leites  integral forms and pseudo-differential forms}
\label{sec:BerLei}

These constructions   were suggested for integration theory on supermanifolds before $r|s$-forms. Their relation with the latter will be considered in the next section.

\subsection{Pseudo-differential forms}
Given a supermanifold $M$ of dimension $n|m$, with local coordinates $x^1, \dots, dx^n$ (even), and $x^{\hat 1}, \dots, x^{\hat m}$ (odd), \emph{pseudo-differential forms} (shortly, \emph{pseudoforms}) on $M$ are  defined  as functions on
$\Pi TM$.  (See Bernstein-Leites~\cite{Bern_Lei77_part2}. Also in~\cite{manin:gaugeeng}, \cite{tv:git}.)
In coordinates, they take the form 
\begin{equation}
	\omega = \omega(x^1, \dots, x^n, x^{\hat 1}, \dots, x^{\hat m}, dx^1, \dots, dx^n, d x^{\hat 1}, \dots, dx^{\hat m}) \, ,
\end{equation}
or, shortly, $\omega=\omega(x,dx)$, where $dx^A$ denote fiber coordinates. The parity of $dx^A$ is opposite to that of $x^A$. The de Rham differential extends to  pseudoforms as an odd vector field $d=dx^A\partial/\partial x^A$ on the supermanifold $\Pi TM$.  

``Ordinary'' differential forms are singled out among pseudoforms by  polynomial dependence on $dx^A$. Dropping  the polynomiality restriction makes it possible to consider objects like $e^{-(dx^{\hat 1})^2}$, i.e. functions on   $\Pi TM$ rapidly decreasing or     compactly-supported along the even directions in the fibers. The point is that such fiberwise rapidly decreasing pseudoforms can be integrated as functions against the Berezin volume element $D(x,dx)$ on   $\Pi TM$,   which can be shown to be invariant. (This is   unlike    differential forms as functions that are fiberwise polynomial.) By definition, this is called the \emph{integral of a pseudoform} over $M$\,:
\begin{equation}\label{eq.intpdf}
  \int_{M}\omega := \int_{\Pi TM}\! D(x,dx)\,\omega(x,dx)\,.
\end{equation}

One observes that in the purely even case $m=0$, all pseudoforms are polynomial is the variables $dx^A$,  which are   odd, and therefore coincide with (inhomogeneous) differential forms. Immediately from the definition of Berezin integral, it follows that in such a case, the integral~\eqref{eq.intpdf} amounts to picking the top-degree part and integrating it in the usual sense as  a top-degree form over $M$.  This is the standard way how integration of inhomogeneous differential forms is defined, hence justifies~\eqref{eq.intpdf} for the general case. 

\subsection{Integral forms} Before introducing pseudo-differential forms, 
Bernstein and Leites  
introduced ``integral forms'' %~\cite{Bern_Lei77_part1}  
in such a way  that they would 
incorporate Berezin volume element $D(x)$, forming a complex dual to the de Rham complex  and allowing to perform coordinate-invariant integration  over supermanifolds and   submanifolds of odd codimension zero, including an analogue of the classical Stokes  theorem. 
In addition to the original papers, see also Witten~\cite{Witten_Notes_on_Super} and \cite{manin:gaugeeng}, \cite{tv:git}.
In physics,  integral forms have been applied to various models, particularly in supergravity and super-Yang–Mills theory (see Grassi et al.~\cite{Grassi2018} and references therein).

The original construction by Bernstein and Leites is as follows.
Let $\Omega^i(M)$ denote the module of polynomial differential $i$-forms on $M$, that is, fiberwise-polynomial functions of degree $i$ from $C^\infty(\Pi TM)$. The  dual module $(\Omega^i(M))^*$ is spanned over $C^\infty(M)$ by polynomials of degree $i$ in partial derivatives with respect to the differentials:
\begin{equation}
x^*_a = \frac{\partial}{\partial dx^a}.
\end{equation}
The space of \emph{integral forms} was defined in~\cite{Bern_Lei77_part1} as the  graded module  
\begin{equation}
\Sigma(M) = \bigoplus_i \Sigma_i(M), \quad \text{where} \quad \Sigma_{n - m - i}(M) = \mathrm{Vol}(M) \otimes_{C^\infty(M)} (\Omega^i(M))^*.
\end{equation}
Hence there is a pairing
\begin{equation*}
  \Omega^i(M)\times \Sigma_{n - m - i}(M) \to \mathrm{Vol}(M)\,.
\end{equation*}
Here  $\mathrm{Vol}(M)$ denotes the space of Berezin  volume forms. There are different conventions for definition of grading for integral forms, see Remark~\ref{rem.gradint}. 
Integral forms  are sums of elements of the form 
\begin{equation} \label{eq:integralform}
	\sigma = f(x) \, D(x) \cdot (x^*_1)^{\alpha_1} \cdots (x^*_n)^{\alpha_n} (x^*_{\hat{1}})^{\beta_1} \cdots (x^*_{\hat{m}})^{\beta_m} \, ,
\end{equation}
where $f(x)$ is a smooth function, $\alpha_i \in \{0,1\}$, and 
$\beta_j \in \mathbb{Z}_{\geq 0}$.

The differential   $d: \Sigma_k \to \Sigma_{k+1}$ is defined by:
\begin{equation}
d  = (-1)^{\tilde a} \frac{\partial}{\partial x^a} \cdot \frac{\partial}{\partial x_a^*}\,,
\end{equation}
so that the action of $\frac{\partial}{\partial x_a^*}$ lowers the power in $x^*_a$, and $\frac{\partial f}{\partial x^a}$ differentiates the coefficient. This differential satisfies $d^2 = 0$ and endows $\Sigma(M)$ with the structure of a cochain complex.

By using fiberwise Fourier transform,   integral forms   can be put into a bijection with the generalized pseudo-differential forms that are  delta-functions in the differentials $dx^A$ and their derivatives, i.e.  fiberwise  distributions on $\Pi TM$ supported at the zero section $M\subset \Pi TM$. Namely,
we have the correspondence:
\[
(x^*_1)^{\alpha_1} \cdots (x^*_n)^{\alpha_n}  (x^*_{\hat{1}})^{\beta_1} \dots  (x^*_{\hat{m}})^{\beta_m} \quad \longleftrightarrow \quad 
(dx^1)^{1 - \alpha_1} \cdots 
(dx^n)^{1 - \alpha_n}
\delta^{(\beta_1)}(dx^{\hat{1}}) \dots 
\delta^{(\beta_m)}(dx^{\hat{m}}) \,, 
\]
up to a factor,  where 
$\alpha_i \in \{0,1\}$ indicates whether the odd variable $dx^i$ is present, and $\delta^{(\beta_j)}(dx^{\hat{j}})$ denotes the $\beta_j$-th derivative of the delta-function of the even variable $dx^{\hat{j}}$.

Altogether, there is the following statement.

\begin{theorem}[Bernstein-Leites~\cite{Bern_Lei77_part2}]
\label{thm.intasgenpdf}
Integral forms on $M$ are isomorphic to  generalized pseudoforms, i.e. generalized functions on $\Pi TM$,
	that are supported on $M\subset \Pi TM$    and  smooth along $M$. This isomorphism is compatible with the differentials and integration. 
\end{theorem}

\begin{remark} \label{rem.gradint}
In the grading above, which was originally proposed in~\cite{Bern_Lei77_part1}, the Berezin volume forms on $M=M^{n|m}$ have degree $n-m$,  
\begin{equation*}
  \mathrm{Vol}(M)=\Sigma_{n-m}(M)\,.
\end{equation*}
This can be motivated by the homogeneity property of the coordinate volume element $D(x,\theta)$  
(here $x$ denotes  the collection of even coordinates and $\theta$,   odd ones) under scaling of the coordinates $x,\theta$. Also, $D(x,\theta)$ corresponds to the generalized pseudoform 
\begin{equation*}
  dx^1 \dots dx^n\, \delta(d\theta^1) \dots \delta(d \theta^m)\,.
\end{equation*}
The latter has homogeneity of weight $n-m$ under scaling transformation of $dx,d\theta$. An alternative grading, which we will denote by upper indices  and adopt in the rest of the paper, is
\begin{equation}
  \Sigma^{i}(M) := \mathrm{Vol}(M) \otimes_{C^\infty(M)} (\Omega^{n-i}(M))^*=\Sigma_{i-m}(M)\,.
\end{equation}
Then the degree assigned to volume forms on $M^{n|m}$ will coincide with   the dimension $n$ of the underlying even manifold,
\begin{equation*}
  \mathrm{Vol}(M)=\Sigma^{n}(M)\,.
\end{equation*}
The elements of $\Sigma^{i}(M)$ are sums of terms~\eqref{eq:integralform} with $i=n-\alpha_1-\ldots-\alpha_n-\beta_1-\ldots\beta_m$. The motivation for this choice is in integration: an integral $n$-form serves for integration over $M^{n|m}$, while an integral $i$-form serves for integration over submanifolds $P^{i|m}\subset M^{n|m}$. %, where $P=P^{i|m}$. 
	(In both definitions of grading, Berezin volume forms are the top forms.)
\end{remark}

An algebraic counterpart of integral forms was introduced in~\cite{Kh_V_recurrent}  as 
\begin{equation}\label{eq.algintf}
  \Sigma^i(V):=\Ber(V)\otimes \Lambda^{n-i}(V^*)\,,
\end{equation}
where $V$ is a vector superspace (we mentioned it in Section~\ref{Sec:E}).  It corresponds to  $\Sigma^i(M)$ when applied to the bundle $T^*M$, up to   choices, e.g. exterior/symmetric powers \,\footnote{Differential forms are seen  in~\cite{Bern_Lei77_part1}  as sections of   $\Lambda^i(T^*M)$ and, in~\cite{Bern_Lei77_part2}, as sections of $S^i(\Pi T^*M)$.}. 

\begin{remark}
\label{rem.orientsec}
  There is a subtlety concerning signs and orientation~\cite{tv:pdf},\cite{tv:git}. %for integrals of pseudo-differential forms and  the relation between integral forms and generalized pseudoforms (see~\cite{tv:pdf} and \cite{tv:git}). 
  The change of variables formula in Berezin integral includes  $\Ber J$ multiplied by $\mathrm{sign}\det J_{00}$ where $J$ is the Jacobi matrix and 
  $J_{00}$ is its even-even block. Therefore, integration of volume forms over a supermanifold $M$ requires   orientation in the sense   $\det J_{00}>0$ between all charts. We may call it \emph{orientation of type I}. By applying this to $\Pi TM$, we deduce that integration of pseudo-differential forms over $M$ requires  for $M$ orientation in a different sense: $\Ber J>0$. We call it \emph{orientation of type II}. To reconcile the two types, one needs a twisting by extra signs on the intersections of charts, namely $\mathrm{sign}\det J_{11}$. So, in Theorem~\ref{thm.intasgenpdf} generalized pseudoforms actually correspond to integral forms with the twisting, which may be called \emph{integral forms of type II}. Respectively, to get the usual integral forms, one has to incorporate the twisting into   the pseudo-differential forms. 
\end{remark}

Since integral forms were designed for integration over $r|m$-dimensional objects  such as $r|m$-paths, it is natural to ask about a relation with $r|m$-forms.   One can obtain from Theorem~\ref{thm:rmforms} an isomorphism of $r|m$-forms and integral $r$-forms for all $r\geq 0$. It  was extended also to   $r<0$, see~\cite{tv:dual}, where $r|s$-forms with negative $r$ were introduced (we do not consider them in this paper).

\section{Baranov-Schwarz transformation}

\label{Sec:BSch}

Baranov and Schwarz~\cite{Baranov_Shwarz} introduced an integral transformation that maps pseudoforms into $r|s$-densities (in their language). 
%Their transformation $\lambda_{r|s}$
%maps Bernstein-Leites integral forms into Vornov-Zorich $r|m$-forms of type II.
% is a way to extract the $r|s$-form component from a general pseudoform $\omega$. It can be viewed as filtering out all but the desired degree from an inhomogeneous expression. 
The following formulation of it was introduced in~\cite{tv:pdf}. Consider an $r|s$-path in $M$ parametrized by   $t \in  U^{r|s}$ for a bounded region $U^{r|s}\subset \mathbb{R}^{r|s}$. Then, by 
expanding $dx(t)$, we obtain $dx(t) = \sum_{F=1}^{r|s} dt^F\, \dot{x}_F^A$, and  define   the \emph{Baranov–Schwarz transformation}  (or \emph{transform})    of a pseudoform $\omega=\omega(x,dx)$,  
\begin{equation}
	\lambda_{r|s} : \omega \mapsto L=\lambda_{r|s} \omega \,,
\end{equation}
where $L=L(x,\dot x)$, by integrating over the variables $dt^F$, as  
\begin{equation}\label{eq.bschint}
	\bigl(\lambda_{r|s} \omega\bigr)(x, \dot{x})=L(x,\dot x) := \int_{\R^{s|r}} D(dt)\, \omega\left(x, \sum_{F=1}^{r|s} dt^F\, \dot{x}_F^A \right) \, .
\end{equation}
Then for every  $r|s$-path  $\Gamma\,:U^{r|s}\to M$, 
\begin{equation}
	\int_{\Gamma} \omega = \int_{U^{r|s}} D(t)\, \int_{\R^{s|r}} D(dt)\, \omega(x(t), dx(t)) = \int_{U^{r|s}} D(t)\, L(x(t), \dot{x}(t)) = \int_{\Gamma} L \,.
\end{equation}
That is, for a pseudoform $\omega$, its integral over an $r|s$-path $\Gamma$ factors into the transformation $\lambda_{r|s}$ and the integral of the Lagrangian $L$. 

\begin{example} On a purely   even manifold, applying the Baranov–Schwarz transformation to an inhomogeneous differential  form returns its homogeneous component.
	More precisely, suppose  $\omega$ is   an inhomogeneous  polynomial in odd variables $dx^i$, where $i \in {1, \dots, n}$.
	Then  $\lambda_{r}(\omega)$ returns the homogeneous component of $\omega$ of degree $r$ considered as a function of $r$ tangent vectors.
 \end{example}

\begin{theorem}[\cite{tv:pdf}] Let $\omega$ be a pseudoform.
	Then the Lagrangian $\lambda_{r|s} \omega$ (when the integral~\eqref{eq.bschint} is defined) is a Voronov-Zorich $r|s$-form of type II.
\end{theorem}

A brief comment about this theorem: in~\cite{tv:pdf} it was proved that Lagrangians $L=\lambda_{r|s}\omega$ satisfy~\eqref{eq.fund}.
In Baranov and Schwarz~\cite{Baranov_Shwarz}, the result of their transformation is already described as an $r|s$-density; however, there is a subtlety concerning the sign. It is of the same nature as mentioned in Remark~\ref{rem.orientsec}. 
Indeed,   a change of parametrization  $t \mapsto t'$ of a path  means substituting $g\dot x$ for $\dot x$ in $L(x,\dot x)$, where $g=\frac{\partial t}{\partial t'}$. One may take arbitrary $g\in GL(r|s)$. %Recall condition~\eqref{eq.bercond} in the definition of Voronov-Zorich forms.
Let us see what will happen  with $L=\lambda_{r|s} \omega$ under such a substitution. We have (suppressing $x$ in the arguments)
\begin{multline*}
  (\lambda_{r|s} \omega)(g\dot x)= \int_{\R^{s|r}} D(dt)\, \omega(\underbrace{dt^Fg_F^H}_{dt'} \dot{x}_H)=
  \int_{\R^{s|r}} D(dt'\cdot g^{-1})\, \omega(dt^{'H} \dot{x}_H)\\
  =\Ber_{1,0}((g^{-1})^{\Pi}) \int_{\R^{s|r}} D(dt')\, \omega(dt^{'H} \dot{x}_H)= \Ber_{0,1}g\cdot (\lambda_{r|s} \omega)(\dot x)\,.
\end{multline*}
Here 
\begin{equation*}
  \Ber_{0,1}g = \operatorname{sign}\det g_{11}\cdot \Ber g\, 
\end{equation*}
and similarly for $\Ber_{1,0}$ (see~\cite{tv:git}). 
This shows that Lagrangians $\lambda_{r|s} \omega$ almost satisfy~\eqref{eq.bercond}, but possibly with a different sign. Such Lagrangians are called \emph{$r|s$-forms of type II} if they satisfy also~\eqref{eq.fund}. They serve for integration over oriented paths with orientation of type II, the same as pseudo-differential forms (see Remark~\ref{rem.orientsec}). 

%with a different notion of orientation defined by the condition $\Ber g>0$ as opposed to $\det g_{00}>0$ needed for integration of Berezin volume forms. 

In the case of $\lambda_{r|m}$ (that is, where the odd degree is maximal), this can be managed with the help of a twist. 
One can show that the transformation $\lambda_{r|m}$, when applied to the pseudoform realizations of Bernstein–Leites integral forms, returns Voronov-Zorich
$r|m$-forms (of type II or usual if   pseudoforms are twisted); and, in fact:   
\begin{theorem}[\cite{tv:git}]
		Bernstein-Leites integral forms of degree $r$ are isomorphic to  Voronov-Zorich
		$r|m$-forms . 
\end{theorem}

In order to illustrate how the Baranov-Schwarz transformation works, consider some elementary examples. %We pick a pseudoform  of a particular appearance. % motivated by the constructions that will be introduced in Section~\ref{sec.geom}.  

Let $M$ be of dimension ${1|1}$ with coordinates $(x, \theta)$,
where $x$ is even and $\theta$ is odd.
Consider a  generalized  pseudoform of the appearance
\[
\omega = f(x, \theta)\cdot \delta(d\theta).
\]
(Recall that $d\theta$ is an even variable.) 

\begin{example}[Baranov--Schwarz transformation $\lambda_{0|1}$]

To compute the Baranov--Schwarz transform $\lambda_{0|1}(\omega)$, we parameterize a  path $\R^{0|1}\to M$ by one odd parameter $\tau$, so the pullback of differentials is: $dx = d\tau \cdot \dot{x}$, $d\theta = d\tau \cdot \dot{\theta}$, 
where $\dot{x}$ is odd and $\dot{\theta}$ is even.
Substitute this into $\omega$:
\[
\omega = f(x, \theta)  \delta(d\tau  \dot{\theta}) = f(x, \theta) \frac{1}{| \dot{\theta} |} \delta(d\tau),
\]
where $|\dot{\theta}|$ is the absolute value. Then:
\[
\lambda_{0|1}(\omega)(x, \theta, \dot{x}, \dot{\theta}) = \int_{\R^{1|0}} D(d\tau)\, f(x, \theta)  \frac{1}{|\dot{\theta}|} \delta(d\tau) =
 \frac{f(x, \theta)}{|\dot{\theta}|} \, .
\]
\end{example}

\begin{example}[Baranov--Schwarz transformation $\lambda_{1|0}$]
Now let the path be given by: $x=x(t)$, $\theta=\theta(t)$, where $t$ is an even parameter. Then $dx=dt\cdot \dot x$, $d\theta = dt \cdot \dot{\theta}$, 
where $\dot x$ is even and $\dot \theta$ is odd. Then,   $\delta(dt \cdot \dot \theta)$ is not defined, and so is $\lambda_{1|0} \omega$.
\end{example}

\begin{example}[Baranov--Schwarz transformation $\lambda_{1|1}$]
	We continue in the settings of the previous examples, but now   $x = x(t, \tau)$, $\theta = \theta(t, \tau)$ with $t$ even, $\tau$ odd. Then $	d\theta = dt \dot{\theta}_t + d\tau \dot{\theta}_\tau$. The Baranov--Schwarz transformation is
	\[
	\lambda_{1|1}(\omega) = \int_{\R^{1|1}} D(dt)\, D(d\tau)\, f(x, \theta) \delta(dt \dot{\theta}_t + d\tau \dot{\theta}_\tau).
	\]
	Using the scaling property of delta-functions in even variables,
	\[
	\delta(d\theta) = \delta(dt \dot{\theta}_t + d\tau  \dot{\theta}_\tau) = \frac{1}{|\dot{\theta}_\tau|} \delta\left(dt \frac{\dot{\theta}_t}{\dot{\theta}_\tau} + d\tau   \right) \, .
	\]
	After a linear change of variables,
	\[
	s = d\tau + dt \frac{\dot{\theta}_t}{\dot{\theta}_\tau}   \quad d t = dt \, , 
	\]
	with the Berezinian equal to $1$,  the integral becomes:
	\[
	\lambda_{1|1}(\omega) = \int_{\R^{1|1}} D(s)\, D(dt)\, 
		\frac{f(x, \theta)}{|\dot{\theta}_\tau|} \delta(s)
	= \int_{\R^{0|1}} D(dt)\, 	
	\frac{f(x, \theta)}{|\dot{\theta}_\tau|} =0 \,  
	\]
because there is no factor of   $dt$ in the integrand. (Alternative calculation can be: $\delta(d\theta) = \delta(dt \dot{\theta}_t + d\tau  \dot{\theta}_\tau) = \frac{1}{|\dot{\theta}_\tau|} \delta\left(dt \frac{\dot{\theta}_t}{\dot{\theta}_\tau} + d\tau   \right)=
\frac{1}{|\dot{\theta}_\tau|} \delta\left(d\tau \right) + 
dt\frac{\dot{\theta}_t}{|\dot{\theta}_\tau|\dot{\theta}_\tau}   \delta'\left( d\tau   \right)
$, and the Berezin integral over $d\tau, dt$ will amount to the integral of $\delta'(d\tau)$ over $d\tau\in \R$, which will give zero.) 
\end{example}

\begin{example}[Baranov--Schwarz transformation $\lambda_{0|2}$] 	We continue in the settings of the previous examples.	Given  two {odd} parameters $\tau_1$ and $\tau_2$, we parameterize $
			\theta(\tau) = \theta_0 + \tau_1 {\theta}_1 + \tau_2 {\theta}_2 + {\tau}_1 {\tau}_2 \theta_{12}$. Then: $d\theta = d\tau_1 \dot{\theta}_1 + d\tau_2 \dot{\theta}_2$, 
			where $\dot{\theta}_1$, $\dot{\theta}_2$ are even. 		
			We want to compute
			\[
			\lambda_{0|2}(f(x, \theta) \delta(d\theta)) = \int D(d\tau_1, d\tau_2) \, f(x, \theta) \delta\bigl(d\tau_1 \dot{\theta}_1 + d\tau_2 \dot{\theta}_2 \bigr).
			\]	
			Using the scaling property of delta functions in even variables,
			\[
			\delta(d\tau_1 \dot{\theta}_1 + d\tau_2 \dot{\theta}_2) = \frac{1}{|\dot{\theta}_1|} \, \delta\left(d\tau_1 + d\tau_2 \frac{\dot{\theta}_2}{\dot{\theta}_1}\right) \, .
			\]
			After a linear change variables,
			\[
			u = d\tau_1 + \frac{\dot{\theta}_2}{\dot{\theta}_1} d\tau_2, \quad v = d\tau_2 \, \
			\]
			with Berezinian equal to $1$.			
			Substituting,
			\[
			\lambda_{0|2}(f(x, \theta) \delta(d\theta)) = \frac{1}{\dot{\theta}_1} \int_{\R^{2|0}} D(u,v) \, f(x, \theta) \delta(u) = \int_{\R^{1|0}} f(x, \theta) dv \, ,
			\]
			which is undefined.
\end{example}

We can observe that in   the above examples, the Baranov-Schwarz transformation applied to the same generalized pseudoform is meaningful only in the case $\lambda_{0|1}$, otherwise   undefined or identically zero\,\footnote{This is by no means a   general situation. One can easily find a rapidly decreasing pseudoform $\omega$, e.g. a bell-shape type, for which $\lambda_{r|s}\omega$ is defined and non-trivial for different ${r|s}$.}. For $\lambda_{0|1}\omega$, the answer can be written as 
\begin{equation*}
  \lambda_{0|1}\omega%(\omega)(x, \theta, \dot{x}, \dot{\theta}) 
                                                         =  \Ber\begin{pmatrix}
                                                                  f & 0 \\
                                                                  \dot x & \dot\theta
                                                                \end{pmatrix}\cdot \mathop{\mathrm{sign}} \dot\theta\,,
\end{equation*}
where 
we have  $(\dot x,\dot\theta)$ as an odd vector in $1|1$-space. % and $\pm=\mathop{\mathrm{sign}} \dot\theta$.
It is an $0|1$-form of type II corresponding to a twisted integral $0$-form on $M^{1|1}$. Compare Theorem~\ref{thm:rmforms}.

\begin{example}[another Baranov--Schwarz transformation $\lambda_{0|1}$]
\label{ex.BS01in11}
If we take 
\begin{equation*}
  \omega= f_1\,\delta(d\theta) + f_{\hat 1}\,dx\,\delta'(d\theta)\,,
\end{equation*}
then we will obtain
\begin{equation*}
  \lambda_{0|1}\omega  =  \Ber\begin{pmatrix}
                                                                  f_1 & f_{\hat 1} \\
                                                                  \dot x & \dot\theta
                                                                \end{pmatrix}\cdot \mathop{\mathrm{sign}} \dot\theta\,,
\end{equation*}
which is a general  $0|1$-form (of type II) in $1|1$-space. It corresponds to an arbitrary integral $0$-form. (In dimension $1|1$, such a form will be  a linear combination of the variables $x^*,\theta^*$.)
\end{example}

In   Section~\ref{sec.conVZ}, we will come back to Baranov--Schwarz transformation and construct its formal analogue for which  there will be no need to consider  forms of type II.

% for particular

\section{More Berezinian expansions}
\label{sec.more}

\subsection{Setup and aim}
%\subsection{Berezinian expansion between the poles}
Let $A$ be an even supermatrix of size $n|m$, considered as an operator on an $n|m$-dimensional superspace $V$. We assume that the entries of $A$ belong to an unspecified Grassmann algebra, so that  it would be better to say that $V$ an $n|m$-dimensional free module over this algebra. We may also treat the entries of $A$ as independent variables.

In Section~\ref{Sec:E}, we recalled the results concerning power expansions of the characteristic function
\[
\Ber(E + zA), 
\]
near $z=0$ and $z=\infty$  (due to Schmitt and Khudaverdian-Voronov).  These can be seen as the two extreme cases of the Laurent expansion at zero. Our goal now will be to study the  Laurent expansions   in the intermediate annular domains, 
with the aim of relating them  to the ``intermediate'' exterior powers $\Lambda^{r|s}(V)$, for $0 < s < m$. 
%Along the way, we will also recover known results due to Schmidt and Khudaverdian--Voronov.

Assume that $A$ is diagonal supermatrix:
\begin{equation}\label{eq.A}
	A=\diag(x_1,\ldots,x_n\,|\,y_1,\ldots,y_m)\,.
\end{equation}
Then: 
\begin{equation}
	 \label{eq.RA}
	\Ber(E+zA)=\frac{\prod_{a=1}^{n}(1+zx_a)}{\prod_{\mu=1}^{m}(1+zy_{\mu})} \, .
\end{equation}
%\begin{multline} \label{eq.RA}
%	\Ber(E+zA)=\Ber\diag(1+zx_1,\ldots,x_1+zx_n\,|\,1+zy_1,\ldots,1+zy_m)=\\
%	\frac{\prod_{a=1}^{n}(1+zx_a)}{\prod_{\mu=1}^{m}(1+zy_{\mu})}= 
%	\prod_{a=1}^{n}(1+zx_a) \prod_{\mu=1}^{m}(1+zy_{\mu})^{-1}\,.
%\end{multline}
The poles are at $z=-1/y_{\mu}$, $\mu=1,\ldots, m$. We assume that all $y_{\mu}$ are invertible and in general position, i.e. all the differences $y_{\mu}-y_{\nu}$ are invertible. To be concrete, assume that
\begin{equation}
	|y_1|<\ldots < |y_m|\,,
\end{equation} 
so, for the poles, we have:
\begin{equation}
	0< \frac{1}{|y_m|}<\ldots < \frac{1}{|y_1|}<\infty\,.
\end{equation}

All inequalities between the elements of a Grassmann algebra are understood as inequalities between the corresponding numerical parts obtained by setting the odd generators to zero. 
The desired Laurent expansions of $\Ber(E + zA)$ are obtained by expanding each factor $(1 + z y_{\mu})^{-1}$ in~\eqref{eq.RA} into a suitable geometric series, and then multiplying all resulting factors through.

\subsection{Formula for expansion near zero}
If $z$ is near zero, i.e. $|z|<\frac{1}{|y_m|}$, that will mean that $|zy_m|<1$ and hence also  $|zy_{\mu}|<1$ for all $\mu$.  
Therefore, we should use the expansion
\begin{equation*}
	(1+zy_{\mu})^{-1}=\sum_{i\geq 0} (-1)^i z^iy_{\mu}^i
\end{equation*}
over the non-negative powers of $z$. This will result in:
\begin{multline*} 
	\Ber(E + zA)=\prod_{a=1}^{n}(1+zx_a) \prod_{\mu=1}^{m}(1+zy_{\mu})^{-1}=
	\prod_{a=1}^{n}(1+zx_a) \prod_{\mu=1}^{m}\left(\sum_{i_{\mu}\geq 0} (-1)^{i_{\mu}} z^{i_{\mu}}y_{\mu}^{i_{\mu}}\right) \, .
\end{multline*}
The first product expands into elementary symmetric polynomials:
\begin{equation*}
	\prod_{a=1}^n (1 + z x_a) = \sum_{k=0}^n z^k \sum_{a_1 < \cdots < a_k} x_{a_1} \cdots x_{a_k} \, .	
\end{equation*}	
The second product is a product of geometric series:
\begin{equation*}
	\prod_{\mu=1}^m \left( \sum_{i_\mu \geq 0} (-1)^{i_\mu} z^{i_\mu} y_\mu^{i_\mu} \right)
	= \sum_{i_1, \ldots, i_m \geq 0} (-1)^{i_1 + \cdots + i_m} z^{i_1 + \cdots + i_m} y_1^{i_1} \cdots y_m^{i_m} \, .
\end{equation*}	
Grouping terms by total degree $\ell = i_1 + \cdots + i_m$ gives:
\begin{equation*}
	\sum_{\ell \geq 0} z^{\ell} (-1)^{\ell} \sum_{\substack{i_1 + \cdots + i_m = \ell}} y_1^{i_1} \cdots y_m^{i_m} \, .
\end{equation*}

So, we have  the Cauchy product of two formal power series:
\begin{multline} \label{eq.zero}	
	\sum_{k=0}^{n}z^k\left(\sum_{a_1<\ldots<a_k} x_{a_1}\ldots x_{a_k}\right)
	\sum_{\ell\geq 0}z^{\ell}(-1)^{\ell}\left(\sum_{\substack{i_1+\ldots+i_m=\ell
			\\i_1\,,\ldots\,, i_m\geq 0}}y_1^{i_1}\ldots y_m^{i_m}\right)=\\
	\sum_{N\geq 0} z^N\left(
	\sum_{\substack{k_1+\ldots+k_n+i_1+\ldots +i_m=N\\ k_1,\,\ldots ,\,k_n=0,1\\
			i_1,\,\ldots , \,i_m\geq 0}}x_1^{k_1}\ldots x_n^{k_n}(-1)^{i_1+\ldots i_m}y_1^{i_1}\ldots y_m^{i_m}\right)\,,
\end{multline}
which involves non-negative powers only. %(We explain the geometrical meaning below.)

\subsection{Formula for expansion near infinity}
At the opposite extreme, when $z$ is near infinity, i.e., $|z| > \frac{1}{|y_1|}$, we have $|z y_{\mu}| > 1$ for all $\mu$, and we should use the expansion
\begin{equation}
	(1 + z y_{\mu})^{-1} = z^{-1} y_{\mu}^{-1} \sum_{j \geq 0} (-1)^j z^{-j} y_{\mu}^{-j}
\end{equation}
for all $\mu$.
 Hence, we obtain:
\begin{multline} \label{eq.inf}
	\Ber(E + zA)=\prod_{a=1}^{n}(1+zx_a) \prod_{\mu=1}^{m}(1+zy_{\mu})^{-1}= 
	\prod_{a=1}^{n}(1+zx_a) \prod_{\mu=1}^{m}\left(z^{-1} y_{\mu}^{-1}\sum_{j_{\mu}\geq 0} (-1)^{j_{\mu}} z^{-j_{\mu}}y_{\mu}^{-j_{\mu}}\right)=\\
	\sum_{k=0}^{n}z^k\left(\sum_{a_1<\ldots<a_k} x_{a_1}\ldots x_{a_k}\right)
	\sum_{p\geq 0}z^{-p-m}(-1)^{p}\left(\sum_{\substack{j_1+\ldots+j_m=p
			\\j_1\,,\ldots\,, j_m\geq 0}}y_1^{-1-j_1}\ldots y_m^{-1-j_m}\right)=\\
	\sum_{N\leq n-m} z^N\left(
	\sum_{\substack{k_1+\ldots+k_n-m-j_1-\ldots -j_m=N \\ k_1,\,\ldots ,\,k_n=0,1\\
			j_1,\,\ldots , \,j_m\geq 0}}x_1^{k_1}\ldots x_n^{k_n}(-1)^{j_1+\ldots+j_m}y_1^{-1-j_1}\ldots y_m^{-1-j_m}\right)\,,
\end{multline} 
where the powers of $z$ are bounded from above by $n-m$. Depending on $n$ and $m$, there can be   negative powers only or also a finite number of positive powers.

These two cases given by formulas~\eqref{eq.zero} and~\eqref{eq.inf} will correspond to the  two cases in Khudaverdian-Voronov theorem (and Schmitt for $s=0$), after we   give them the geometric interpretations.  

\subsection{Formula for ``intermediate'' expansion}
Now let us consider the general case of the expansion of $\Ber(E + zA)$ in the annulus between two poles, which was not previously considered. Consider some $s$, where $0<s<m$. Let 
\begin{equation*}
	\frac{1}{|y_{m-s+1}|}< |z| < \frac{1}{|y_{m-s}|}\,.
\end{equation*}
That means that $|zy_1|<1\,,\ldots, |zy_{m-s}|<1$, and
$|zy_{m-s+1}|>1\,,\ldots,|zy_{m}|>1$. Therefore, we  should use the expansion
\begin{equation}
	(1+zy_{\mu})^{-1}=\sum_{i\geq 0} (-1)^i z^iy_{\mu}^i
\end{equation}
for $\mu=1,\ldots,m-s$ and the expansion
\begin{equation}
	(1+zy_{\mu})^{-1}=z^{-1} y_{\mu}^{-1}\sum_{j\geq 0} (-1)^j z^{-j}y_{\mu}^{-j}
\end{equation}
for $\mu=m-s+1,\ldots,m$\,. By substituting into~\eqref{eq.RA}, we obtain:
\begin{multline} \label{eq.med}
	\Ber(E + zA)=\prod_{a=1}^{n}(1+zx_a) \prod_{\mu=1}^{m}(1+zy_{\mu})^{-1}\\
	=\prod_{a=1}^{n}(1+zx_a) \prod_{\mu=1}^{m-s}\left(\sum_{i_{\mu}\geq 0} (-1)^{i_{\mu}} z^{i_{\mu}}y_{\mu}^{i_{\mu}}\right)
	\prod_{\nu=m-s+1}^{m}\left(z^{-1} y_{\nu}^{-1}\sum_{j_{\nu}\geq 0} (-1)^{j_{\nu}} z^{-j_{\nu}}y_{\nu}^{-j_{\nu}}\right)\\
	=\sum_{k=0}^{n}z^k\left(\sum_{a_1<\ldots<a_k} x_{a_1}\ldots x_{a_k}\right)
	\sum_{\ell\geq 0}z^{\ell}(-1)^{\ell}\left(\sum_{\substack{i_1+\ldots+i_{m-s}=\ell
			\\i_1\,,\ldots\,, i_{m-s}\geq 0}}y_1^{i_1}\ldots y_{m-s}^{i_{m-s}}\right)\\
	\times \sum_{p\geq 0}z^{-p-s}(-1)^{p}\left(\sum_{\substack{j_{m-s+1}+\ldots+j_m=p
			\\j_{m-s+1}\,,\ldots\,, j_m\geq 0}}y_{m-s+1}^{-1-j_{m-s+1}}\ldots y_m^{-1-j_m}\right)\\
	=\sum_{N\in \ZZ} z^N\left(
	\sum_{\substack{k+\ell-p-s=N \\ 
			k_1+\ldots+k_n=k\\
			i_1+\ldots+i_{m-s}=\ell\\
			j_{m-s+1}+\ldots+j_m=p\\
			k_1,\,\ldots ,\,k_n=0,1\\
			i_1\,,\ldots\,, i_{m-s}\geq 0\\
			j_1,\,\ldots , \,j_m\geq 0}} 
	x_1^{k_1}\ldots x_n^{k_n}\,(-1)^{\ell+p}\, y_1^{i_1}\ldots y_{m-s}^{i_{m-s}}\, y_{m-s+1}^{-1-j_{m-s+1}}\ldots y_m^{-1-j_m}\right)\,.
\end{multline}

As we will see, for the sake of the geometric interpretation that follows, it is convenient to include the sign factor $(-1)^N$ in the coefficient of the $N$th power of $z$, and thus rewrite the formula as
\begin{multline} \label{eq.med2}
	\Ber(E + zA)=\\
	\sum_{N\in \ZZ} z^N(-1)^N
	\left(
	\sum_{\substack{k+\ell-p-s=N \\ 
			k_1+\ldots+k_n=k\\
			i_1+\ldots+i_{m-s}=\ell\\
			j_{m-s+1}+\ldots+j_m=p\\
			k_1,\,\ldots ,\,k_n=0,1\\
			i_1\,,\ldots\,, i_{m-s}\geq 0\\
			j_1,\,\ldots , \,j_m\geq 0}} 
	(-1)^{k+s}\,x_1^{k_1}\ldots x_n^{k_n}\,y_1^{i_1}\ldots y_{m-s}^{i_{m-s}}\,\,y_{m-s+1}^{-1-j_{m-s+1}}\ldots y_m^{-1-j_m}\right)\,.
\end{multline}

In particular, for the special cases of~\eqref{eq.zero} and \eqref{eq.inf}, we will have, respectively,
\begin{equation}\label{eq.nearzero}
  \Ber(E + zA)=\sum_{N\geq 0} z^N(-1)^N
  \left(
	\sum_{\substack{k+\ell=N\\
k_1+\ldots+k_n=k\\
i_1+\ldots+i_{m}=\ell\\ 
k_1,\,\ldots ,\,k_n=0,1\\
			i_1,\,\ldots , \,i_m\geq 0}}
(-1)^{k}x_1^{k_1}\ldots x_n^{k_n}y_1^{i_1}\ldots y_m^{i_m}\right)\,,
\end{equation}
near zero, and
\begin{equation}\label{eq.nearinfin}
  \Ber(E + zA)=\sum_{N\leq n-m} z^N(-1)^N
  \left(
	\sum_{\substack{k-p-m=N \\ 
k_1+\ldots+k_n=k\\
j_{1}+\ldots+j_m=p\\
k_1,\,\ldots ,\,k_n=0,1\\
			j_1,\,\ldots , \,j_m\geq 0}}
(-1)^{k+m}\,x_1^{k_1}\ldots x_n^{k_n}y_1^{-1-j_1}\ldots y_m^{-1-j_m}\right)\,,
\end{equation}
near infinity.

\section{Geometric interpretation} \label{sec.geom}

For our vector space $V$, consider the parity-reversed space $\Pi V$. The motivation is to work with symmetric powers instead of exterior powers, 
because symmetric powers are easier to handle since their elements can be seen as functions. The said holds for the known case of $\Lambda^k(V)$, where we have $\Lambda^k(V)\cong \Pi^kS^k(\Pi V)$, and $S^k(\Pi V)$ is a subspace of functions on $\Pi V^*$ regarded as a supermanifold; but this is also our plan for the ``intermediate'' exterior powers $\Lambda^{r|s}(V)$, which we want to relate with the   expansions of the characteristic function.    

If a basis in $V$ is denoted by $e_A$, more specifically $e_a, e_{\mu}$, where $a = 1, \ldots, n$ and $\mu = \hat{1}, \ldots, \hat{m}$, then the corresponding basis in $\Pi V$ will be denoted by $\e_A$, where $\e_A = \Pi e_A$. Hence in $\Pi V$ we will have odd basis vectors $\e_1, \ldots, \e_n$ and even basis vectors $\e_{\hat{1}}, \ldots, \e_{\hat{m}}$.

Let us assume that the basis $e_A \in V$, and hence the basis $\e_A \in \Pi V$, consists of eigenvectors of the operator $A$. (Recall that we assume $A$ to be even, invertible, and in general position in the sense explained above.) We have:
\begin{equation}
	A(\e_a)= x_a\e_a\,\quad \text{and} \quad A(\e_{\mu})=y_{\mu} \e_{\mu}\,,
\end{equation}
for some even eigenvalues $x_a$, $y_{\mu}$.

\subsection{Expansion near zero.}

Consider the symmetric powers $S^r(\Pi V)$\,,  $r=0,1,2,\ldots\ $. From the basis $\e_A\in \Pi V$, we obtain the following basis in the space $S^r(\Pi V)$\,: 
\begin{equation}
	(\e_1)^{k_1}\ldots (\e_n)^{k_n}(\e_{\hat 1})^{i_1}\ldots (\e_{\hat m})^{i_m}\,,\label{eq.basinsr}
\end{equation}
where $k_a\in\{0,1\}$\,, $i_{\mu}\in\ZZ_{\geq 0}$\,, and $k_1+\ldots k_n+i_{1}+\ldots i_{m}=r$. The parity of the element~\eqref{eq.basinsr} equals $k=k_1+\ldots k_n \mod 2$. 

For example, in $S^2(\Pi V)$, we obtain the basis consisting  of the elements $\e_a\e_b$, $a<b$  (even), $\e_a\e_{\mu}$ (odd), and $\e_{\mu}\e_{\nu}$, $\mu\leq \nu$ (even). 

If we consider the natural action of $A$ induced on $S^r(\Pi V)$\,---\,that is, more formally, the symmetric power of the parity-reversed operator $S^r(A^{\Pi})$\,---\,then on the basis vectors~\eqref{eq.basinsr}  we will obtain
\begin{multline} \label{eq.srbasis}
	S^r(A^{\Pi})\left((\e_1)^{k_1}\ldots (\e_n)^{k_n}(\e_{\hat 1})^{i_1}\ldots (\e_{\hat m})^{i_m}\right)= \\
	x_1^{k_1}\ldots x_n^{k_n}y_{1}^{i_1}\ldots y_{m}^{i_m}\cdot 
	(\e_1)^{k_1}\ldots (\e_n)^{k_n}(\e_{\hat 1})^{i_1}\ldots (\e_{\hat m})^{i_m}\,,
\end{multline}
i.e., we have eigenvectors with the eigenvalues 
\begin{equation}
	x_1^{k_1}\ldots x_n^{k_n}y_{1}^{i_1}\ldots y_{m}^{i_m}\,.
\end{equation}
Hence,
\begin{equation}
	\str S^r(A^{\Pi})= \sum_{k+\ell=r} (-1)^k\,\sum_{\substack{k_1+\ldots k_n=k\\
			i_{1}+\ldots i_{m}=\ell}}x_1^{k_1}\ldots x_n^{k_n}y_{1}^{i_1}\ldots y_{m}^{i_m}\,.
\end{equation}
Comparing the answer with   \eqref{eq.zero} and \eqref{eq.med2}, we arrive at the following theorem.
\begin{theorem}
	\begin{equation}
		\Ber(E + zA)=\sum_{r\geq 0} z^r(-1)^r\str S^r(A^{\Pi})
	\end{equation}
	for $z$ near zero.
	%\qed
\end{theorem}

From the isomorphism $S^r(\Pi V)\cong \Pi^r\Lambda^r(V)$,  we have that $S^r(A^{\Pi})$ is similar to $\Lambda^r(A)$ for even $r$ and 
to $\left(\Lambda^r(A)\right)^{\Pi}$ for odd $r$. 
Since, for any $B$,  $\str(B^{\Pi})=-\str B$, we see that the above statement gives the theorem of Schmitt and Khudaverdian-Voronov on the expansion of $\Ber(E + zA)$ at zero.

\subsection{Expansions near infinity and intermediate.}

To obtain the geometric interpretation of the expansions near infinity and in the intermediate annuli, we introduce the language of \textit{formal delta-functions}. Note: the term ``formal delta-function'' may refer to different notions; the one used here seems to be not the same as   in vertex algebra theory.

\subsubsection{Formal delta-functions.}\label{sec.fdeltas}
Let $t$ be an even variable, which we may regard as complex. The symbol $\delta(t)$ is defined by the following axiomatic properties:
\begin{align}
	&t\,\delta(t) =0 \label{eq.deltavac}\\
	&	\delta(a\,t) =\frac{1}{a}\,\delta(t) \quad \text{for any even invertible $a$}\label{eq.deltahom}\\
	&	\text{$\delta(t)$ is odd}\,.\label{eq.deltaodd}
\end{align} 

We also introduce the symbols $\delta'(t), \delta''(t), \dots, \delta^{(k)}(t)$ (all odd), interpreted as the derivatives of the formal delta-function $\delta(t)$, and assume that the usual rules of differentiation apply.  For example, we obtain, for the product with $t$,
\begin{equation}
	\delta(t)+t\,\delta'(t)=0\,,
\end{equation}
and we also obtain that
\begin{equation}
	\delta'(a\,t) =\frac{1}{a^2}\,\delta(t)\,.
\end{equation}
By induction,
\begin{equation}
	\delta^{(k)}(a\,t) =\frac{1}{a^{k+1}}\,\delta^{(k)}(t)\,.\label{eq.deltakhom}
\end{equation}
We can introduce these formal symbols and take their linear combinations. We can also compute products of these formal delta-functions and their derivatives, provided they have different arguments --- for example, $\delta(t_1)\,\delta'(t_2)$ for independent variables $t_1$ and $t_2$.

One can see that the postulated properties~\eqref{eq.deltavac} and \eqref{eq.deltahom} of these formal symbols $\delta(t)$ resemble those of the actual Dirac delta-function. The key difference is that for the actual delta-function, $\delta(at) = \frac{1}{|a|},\delta(t)$, whereas in~\eqref{eq.deltahom}, the absolute value is absent.

Formal delta-functions with the above properties were considered by Belopolsky~\cite{belopolsky1997new} and Witten~\cite{Witten_Notes_on_Super}, who also argued that it is natural to treat them as odd quantities.

Properties~\eqref{eq.deltavac}, \eqref{eq.deltahom}, and \eqref{eq.deltaodd} admit the following \textbf{interpretation}: to an even variable $t$, we associate a Fock space, where $\delta(t)$ plays the role of the vacuum vector, and multiplication by $t$ and differentiation with respect to $t$ correspond to the annihilation and creation operators, respectively. We may consider ``finite'' vectors in this Fock space (i.e., finite linear combinations of $\delta(t)$ and its derivatives), as well as ``infinite'' vectors (infinite formal series). One of the questions we aim to address is the behavior of this Fock space\,---\,and tensor products of such spaces\,---\,under substitutions in the variable(s) $t$ (or $t_1, t_2, \ldots$). In this picture, $\delta(t)$ does not come about as odd automatically, but we can apply the parity reversion functor. 

\begin{remark}
\label{rem.remodess}
  A realization of formal delta-functions was suggested to us by A.~Odesskii (private communication). Consider the space of formal power series $\C[[t,t^{-1}]]$ as a module over the polynomial algebra $\C[t]$. Consider the quotient module $M_t:=\C[[t,t^{-1}]]/\C[[t]]$. Then the element $[t^{-1}]\in M_t$ satisfies properties~\eqref{eq.deltavac} and \eqref{eq.deltahom}. Denote  the usual derivative in $t$ by $\partial$;  then $\C[[t]]\subset \C[[t,t^{-1}]]$ is a differential submodule with $\partial$ as the derivation, and $M_t$ will inherit the differential module structure. 
  Ultimately, we  obtain $M_t$ as our Fock space (apart from parity reversion)  with  the class $[t^{-1}]$ as the vacuum vector, $t$ as the annihilation operator and $\partial$ as the creation operator.  So, here  $\delta(t)$ will be represented by   $[t^{-1}]$   and $\delta^{(k)}(t)$, for all $k\in\Z_{>0}$,   by $\partial^k[t^{-1}]=(-1)^kk![t^{-1-k}]$. To have also~\eqref{eq.deltaodd}, we can additionally multiply by $\Pi$ from the left. 
\end{remark}

The usefulness of the above interpretation is in that it also allows to introduce formal integration rules (see later in Section~\ref{sec.conVZ}).

\subsubsection{Vector space corresponding to expansion near infinity.}

%We   continue to use the   notation $\e_a$ (odd vectors) and $\e_{\mu}$ (even vectors) for a  basis in $\Pi V$, which is assumed to be an eigenbasis of the operator $A$. (Everything as before.) 

From a basis $\e_A$ of $\Pi V$, 
define a new vector space with a basis
\begin{equation}
	(\e_1)^{k_1}\ldots (\e_n)^{k_n}
	\delta^{(j_1)}(\e_{\hat 1})\ldots \delta^{(j_m)}(\e_{\hat m})\,,
	\label{eq.basinsigma}
\end{equation} 
so that
\begin{equation}\label{eq.weight}
	k_1+\ldots+k_n-m-j_1-\ldots -j_m=N\,,
\end{equation}
for a fixed number $N\in \ZZ$, 
where the deltas are the formal delta-functions defined above. Here, the $\e_{\mu}$ are regarded as independent even variables. The parity of the basis vector~\eqref{eq.basinsigma} is given by
\begin{equation}\label{eq.pari}
	k_1+\ldots+k_n+m\, \mod 2 \, .
\end{equation}
%since each delta-function is odd.
We denote the vector space defined above by $S^{N,m}(\Pi V)$. From condition~\eqref{eq.weight}, it follows that $N\leq n-m$.  We set $S^{N,m}(\Pi V)=0$ for $N>n-m$.

Under the  scaling $\e_A\mapsto \lambda e_A$, where $\lambda$ is an even invertible scalar, the basis vectors~\eqref{eq.basinsigma} transform as follows:
\begin{equation}\label{eq.scal}
	(\e_1)^{k_1}\ldots (\e_n)^{k_n}
	\delta^{(j_1)}(\e_{\hat 1})
	\ldots 
	\delta^{(j_m)}(\e_{\hat m}) \mapsto \lambda^N\, (\e_1)^{k_1}\ldots (\e_n)^{k_n}
	\delta^{(j_1)}(\e_{\hat 1})
	\ldots 
	\delta^{(j_m)}(\e_{\hat m})\, 
\end{equation}
as a consequence of conditions~\eqref{eq.weight} and \eqref{eq.deltakhom}. 
Thus, \emph{all elements of the vector  space $S^{N,m}(\Pi V)$  are homogeneous of   weight  $N$  with respect to   scaling in $V$ or $\Pi V$. }

The   action of $A$ on $V$ (and the induced action as $A^{\Pi}$ on $\Pi V$) induces the action on the basis vectors~\eqref{eq.basinsigma} and, consequently, to the space $S^{N,m}(\Pi V)$, which we denote by  $S^{N,m}(A^{\Pi})$\,:
\begin{multline} \label{eq.sNmA}
	S^{N,m}(A^{\Pi})\left((\e_1)^{k_1}\ldots (\e_n)^{k_n}
	\delta^{(j_1)}(\e_{\hat 1})
	\ldots 
	\delta^{(j_m)}(\e_{\hat m})\right)= \\
	x_1^{k_1}\ldots x_n^{k_n}y_1^{-1-j_1}\ldots y_m^{-1-j_m}\cdot 
	(\e_1)^{k_1}\ldots (\e_n)^{k_n}
	\delta^{(j_1)}(\e_{\hat 1})
	\ldots 
	\delta^{(j_m)}(\e_{\hat m})\,,
	\,,
\end{multline}
i.e., they will be eigenvectors with the eigenvalues 
\begin{equation}
	x_1^{k_1}\ldots x_n^{k_n}y_1^{-1-j_1}\ldots y_m^{-1-j_m}\,.
\end{equation}
Hence:
\begin{equation}
	\str S^{N,m}(A^{\Pi})= \sum_{k+ p-m=N} 
	\sum_{\substack{k_1+\ldots+k_n=k\\
			j_{1}+\ldots+j_m=p}} 
	(-1)^{k+m}\,x_1^{k_1}\ldots x_n^{k_n}y_1^{-1-j_1}\ldots y_m^{-1-j_m}\, ,
\end{equation}
where $k_1,\,\ldots ,\,k_n=0,1$ and 	$j_1,\,\ldots , \,j_m\geq 0$. By comparing with~\eqref{eq.nearinfin}, we arrive at the following statement:
%with~\eqref{eq.inf} and observing that $(-1)^{j_1+\ldots+j_m} = (-1)^N (-1)^{k+m}$, we arrive at the following statement:
\begin{theorem}
	For $z$ near infinity,
	\begin{equation}\label{eq.expinf}
		\Ber(E + zA)=\sum_{N\leq n-m} z^N(-1)^N\str S^{N,m}(A^{\Pi})\,. %\qed
	\end{equation}
\end{theorem}
Here the convention regarding the parity of the formal delta-function plays   crucial role in determining the signs in the supertrace.

To compare with the Khudaverdian-Voronov result, we need to compare the vector spaces $S^{N,m}(\Pi V)$ with the spaces $\Ber V \otimes \Lambda^r(V^*)$. In brief, they are isomorphic up to parity reversion when $N = n - m - r$. This relation between $N$ and $r$ follows from comparing the weights: $N$ for $S^{N,m}(\Pi V)$ and $n - m - r$ for $\Ber V \otimes \Lambda^r(V^*)$.

This isomorphism can be obtained by a ``formal Fourier transform'' between functions of the variables $\e_A$ and functions of the dual variables $\e^A$ (the basis vectors in $\Pi V^*$). If we assume that Fourier transform acts on formal delta-functions in the same way as on ordinary delta-functions, then  it will map the vector~\eqref{eq.basinsigma} to the monomial
\begin{equation}
	(\e^1)^{1-k_1}\ldots (\e^n)^{1-k_n}
	(\e_{\hat 1})^{j_1}
	\ldots 
	(\e_{\hat m})^{j_m}\in S^{n-k+p}(\Pi V^*)
	\label{eq.basinsigma2}
\end{equation} 
(up to a factor), where $k_1+\ldots+k_n=k$ and  $j_{1}+\ldots+j_m=p$. Since $k-m-p=N$, we have $n-k+p=n-m-N$, so we are in $S^{n-m-N}(\Pi V^*)$\,
\footnote{For the Fourier transform of functions on a vector space to be coordinate-independent, the image has to be balanced by the factor of the coordinate volume element, which is in our case $D(\e^{\boldsymbol{\cdot}})\in \Vol(\Pi V)\cong \Ber \Pi V^*\cong \Ber V$. Hence, more precisely, we get an element of $S^{n-m-N}(\Pi V^*) \otimes \Ber V\,$.}. 

We need to take into account that the formal Fourier transform in question has parity $n-m \mod 2$ (because there are $n$ odd variables of integration and   there is a product of $m$ odd formal delta-functions). Hence we arrive at an  even  isomorphism
\begin{equation*}
	S^{N,m}(\Pi V)\cong \Pi^{n-m}\,S^{n-m-N}(\Pi V^*) \otimes \Ber V\,.
\end{equation*}
By composing it with the isomorphism  $S^{n-m-N}(\Pi V^*)\cong \Pi^{n-m-N} \Lambda^{n-m-N}(V^*)$, we  obtain 
\begin{equation} \label{eq.iso}
	S^{N,m}(\Pi V)\cong \Pi^{N}\,\Lambda^{n-m-N}(V^*) \otimes \Ber V\,,
\end{equation}
which is the exact claim. Compare with~\eqref{eq.algintf} and \eqref{exp_infty}.

\begin{remark} In the context of forms on supermanifolds and in the   language of actual distributions instead of formal delta-functions used here, this corresponds to Theorem~\ref{thm.intasgenpdf}. %it is a well-known result about the isomorphism of Bernstein-Leites integral forms on $M$ with the generalized pseudo-differential form supported on $M\subset \Pi TM$\,. (See Bernstein-Leites~\cite{Bern_Lei77_part2} and Voronov~\cite{tv:git}.)
\end{remark}

\begin{remark}\label{rem.substdeltas}
	If we consider general  non-degenerate  linear substitutions into formal delta-functions, we can observe that, for the product $\delta(t_1)\ldots \delta(t_m)$,
	\begin{equation}
		\delta(t_{\mu_1}T^{\mu_1}_1)\ldots \delta(t_{\mu_m}T^{\mu_m}_m)=
		\delta(t_1)\ldots \delta(t_m)\,(\det T)^{-1}\,.
	\end{equation}
	Indeed, one can consider $T$ as a product of diagonal and elementary matrices. For diagonal matrices, we use~\eqref{eq.deltahom}; and for an elementary matrix, which gives a substitution  of the form $\delta(t_{\mu}+c t_{\nu})$ into one delta-function $\delta(t_{\mu})$ and keeps the rest unchanged, we see that it will leave  the product unchanged due to~\eqref{eq.deltavac} (consider the Taylor expansion at $t_{\mu}$ and observe that all the extra terms will be annihilated by multiplication with $\delta(t_{\nu})$). More generally, if we include a set of odd variables $\theta_1,\ldots,\theta_n$ and consider the product
	\begin{equation*}
		\theta_1\ldots \theta_n \delta(t_1)\ldots \delta(t_m)\,,
	\end{equation*}
	which is the formal delta-function of the whole set of $t$'s and $\theta$'s\,\footnote{For an odd variable, the formal delta-function is the same as the usual delta-function.},
	we will get a Berezinian in the denominator by the same argument. Specifically, in our geometric case, if we denote
	\begin{equation*}
		\delta(\e_{\boldsymbol{.}}):=\e_1\ldots \e_n \delta(\e_{\hat 1})\ldots \delta(\e_{\hat m})\,,
	\end{equation*}
	the formal delta-function of the   variables $\e_A=(\e_{\mu},\e_a)$, 
	and consider a change of basis $e_{A'}=e_{A}T^A_{A'}$ in $V$, which induces $\e_{A'}=\e_{A}T^A_{A'}$ in $\Pi V$, then
	\begin{equation}
		\delta(\e_{\boldsymbol{.}}T) = \delta(\e'_{\boldsymbol{.}}) \Ber T\, 
	\end{equation}
	($\Ber T$ in the numerator because of the parity reversion). This shows that the formal delta-function $\delta(\e_{\boldsymbol{.}})$ can be identified with the basis element $\Ber(e)\in \Ber V$ corresponding to the basis $e_A\in V$. Following this way, one can obtain  another proof of the isomorphism~\eqref{eq.iso}.
	
\end{remark}

\subsubsection{Vector space corresponding to ``intermediate'' expansion.}

Fix an $s$ with $0 < s < m$. Continuing with the same notation as above, we introduce a new vector space $S^{N,s}(\Pi V)$ spanned by the products 
\begin{equation}
	(\e_1)^{k_1}\ldots (\e_n)^{k_n}
	(\e_{\hat 1})^{i_1}\ldots (\e_{\widehat {m-s}})^{i_{{m-s}}}
	\delta^{(j_{m-s+1})}(\e_{\widehat{m-s+1}})\ldots \delta^{(j_m)}(\e_{\hat m})\,, \label{eq.basinsNs}
\end{equation}
where 
\begin{equation}\label{eq.weights}
	k_1+\ldots+k_n
	+i_1+\ldots+i_{m-s} 
	-s-j_{m-s+1}-\ldots -j_m=N\,,
\end{equation}
as a basis, for a fixed number $N\in \ZZ$. Now $N$ can be any integer without restrictions. %, positive, zero, or negative. 

Similarly to the above, we observe that all elements of $S^{N,s}(\Pi V)$ have weight $N\in \ZZ$ with respect to scaling transformations of $V$. The parity of the basis vector~\eqref{eq.basinsNs} is 
\begin{equation}\label{eq.paris}
	k_1+\ldots+k_n+s\, \mod 2\,.
\end{equation}

If we allow  $s=0$ and $s=m$, the construction recovers the previous case. In particular,  $S^{N,0}(\Pi V)=S^N(\Pi V)$. Unlike these two cases of $s=0$ or $s=m$, the definition of the vector space $S^{N,s}(\Pi V)$ for $0<s<m$ explicitly depends on the ordering of a basis $e_A\in V$. Hence there is the problem of its basis (in)dependence, which we do not fully address here, but will consider in detail elsewhere.

The   action of $A$ on $V$ (and hence on $\Pi V$) induces the action on the space $S^{N,s}(\Pi V)$,  which we denote    $S^{N,s}(A^{\Pi})$. 
This is defined by its induced action
on the basis vectors~\eqref{eq.basinsNs}\,:
\begin{multline} \label{eq.sNsA}
	S^{N,s}(A^{\Pi})\,: (\e_1)^{k_1}\ldots (\e_n)^{k_n}
	(\e_{\hat 1})^{i_1}\ldots (\e_{\widehat {m-s}})^{i_{{m-s}}}
	\delta^{(j_{m-s+1})}(\e_{\widehat{m-s+1}})\ldots \delta^{(j_m)}(\e_{\hat m}) \mapsto \\
	x_1^{k_1}\ldots x_n^{k_n}\,y_1^{i_1}\ldots y_{m-s}^{i_{m-s}}\, y_{m-s+1}^{-1-j_{m-s+1}}\ldots y_m^{-1-j_m} {\quad \quad \quad \quad \quad \quad \quad \quad \quad \quad }
	\\ \cdot  
	(\e_1)^{k_1}\ldots (\e_n)^{k_n}
	(\e_{\hat 1})^{i_1}\ldots (\e_{\widehat {m-s}})^{i_{{m-s}}}
	\delta^{(j_{m-s+1})}(\e_{\widehat{m-s+1}})\ldots \delta^{(j_m)}(\e_{\hat m})\,,
\end{multline}
which will be eigenvectors with the eigenvalues 
\begin{equation}
	x_1^{k_1}\ldots x_n^{k_n}\,y_1^{i_1}\ldots y_{m-s}^{i_{m-s}}\,y_{m-s+1}^{-1-j_{m-s+1}}\ldots y_m^{-1-j_m}\,.
\end{equation}

Hence
\begin{equation}
	\str S^{N,s}(A^{\Pi})= \sum_{k+\ell-p-s=N} 
	(-1)^{k+s} 
	\sum_{\substack{k_1+\ldots+k_n=k\\
			i_1+\ldots+i_{m-s}=\ell\\
			j_{m-s+1}+\ldots+j_m=p}} 
	x_1^{k_1}\ldots x_n^{k_n}\,y_1^{i_1}\ldots y_{m-s}^{i_{m-s}}\,y_{m-s+1}^{-1-j_{m-s+1}}\ldots y_m^{-1-j_m}\,,
\end{equation}
with $k_a\in\{0,1\}$ and $i_{\mu}, j_{\nu}\geq 0$. By comparing with~\eqref{eq.med2}, we obtain the following statement:
\begin{theorem}
\label{thm.geomint}
	For any $s$, $0<s<m$, the Laurent expansion of $\Ber(E + zA)$ in the annulus 
	\begin{equation}
		\frac{1}{|y_{m-s+1}|}< |z| < \frac{1}{|y_{m-s}|} 
	\end{equation}
	between two consecutive poles  
	is given by
	\begin{equation}\label{eq.medgeom}
		\Ber(E + zA)=\sum_{N\in \ZZ}z^N(-1)^N\str S^{N,s}(A^{\Pi}) = \sum_{N\in \ZZ} z^N\str \left(S^{N,s}(A^{\Pi})\right)^{\Pi^N}\,. 
	\end{equation}
	The coefficients in the right-hand side are the supertraces of the action in the vector spaces $\Pi^N S^{N,s}(\Pi V)$. \qed
\end{theorem}

\begin{remark}
	We have obtained the theorem using a construction of the space $S^{N,s}(\Pi V)$ (a ``generalized'' symmetric power of $\Pi V$, which should correspond to a ``generalized'' exterior power $\Lambda^{r|s}(V)$) with the help of a particular basis in $V$ in which the operator $A$ is diagonal. Hence we need to address the problem of the basis dependence as well as defining the action of a non-diagonal $A$ on vectors of the form~\eqref{eq.basinsNs}. These are two sides of the same problem, and, in brief, a solution can be obtained by considering $S^{N,s}(\Pi V)$ as a Fock space, see~\ref{sec.fdeltas}, and constructing there a   Berezin-style spinor representation\,\footnote{As Th.~Voronov   has  suggested to us, a choice of splitting of the   basis vectors $e_{\mu}$ into   groups of $m-s$ and $s$ elements is similar to a polarization in geometric quantization (private communication).}. We will consider it elsewhere.
\end{remark}

\section{Connection with Voronov-Zorich forms.} \label{sec.conVZ}

%\subsection{Formal Baranov-Schwarz transformation and connection with Voronov-Zorich forms.}

The vector spaces $S^{N,s}(\Pi V)$ introduced above can be   connected with Voronov-Zorich  forms as follows. In order to avoid re-stating definitions from Section~\ref{Sec:VZ} in the  setup for multivectors (this can be found in~\cite{shemya:voronov2019:superplucker}, see also appendix to the second edition of~\cite{tv:git}), we change the notation here and assume that the basis of the vector superspace $V$ consists of the differentials of coordinates on a supermanifold $M$ and will therefore write it with upper instead of lower indices. (Thus we will be speaking about ``forms'' instead of ``multivectors'' and employ differential-geometric notation.) Hence the formal delta-functions will be  %$\delta(dx^a)$ and 
$\delta(d\xi^{\mu})$, where $x^a$ and $\xi^{\mu}$ are even and odd coordinates, respectively. Since $dx^a$ are odd and  $\delta(dx^a)=dx^a$,   we may also speak about formal delta-functions of both $dx^a, d\xi^{\mu}$ for uniformity.  The expression of the form
\begin{equation}\label{eq.fpdf}
	(dx^1)^{k_1}\,\ldots\, (dx^n)^{k_n}\,
	(d\xi^1)^{i_1}\ldots (d\xi^{m-s})^{i_{{m-s}}}\,
	\delta^{(j_1)}(d\xi^{m-s+1})\,\ldots \,\delta^{(j_s)}(d\xi^m)  %\label{eq.basinsNs}
\end{equation}
can be regarded as a formal (generalized) pseudo-differential form. (Compare considerations in Sections~\ref{sec:BerLei} and~\ref{Sec:BSch}.) Its definition depends on a splitting of the odd coordinates into two groups, the first $m-s$ and the last $s$. 

\subsection{Formal integration.}
Our aim is to introduce a formal analogue of the Baranov-Schwarz transformation mimicking the construction considered in Section~\ref{Sec:BSch}. To this end, we need to define integration of the expressions~\eqref{eq.fpdf} and those obtained from them by substitutions. This will be called ``formal integration''. In a nutshell, we can do that by postulating that the integral of the formal delta-function is calculated by the same rule as for the ordinary delta-function (except for a change of variable),
\begin{equation}\label{eq.intdelta}
  \fint  \delta(t)\,dt=1\\, \quad \fint  \delta'(t)\,dt=0\,, \quad \text{etc.}
\end{equation} 
(We   use  a modified integral symbol and  do not indicate a ``region of integration''.)

In view of realization of formal delta-functions offered by Remark~\ref{rem.remodess}, this is more than a postulate. 
Since, for a single variable $t$, the space $M_t=\C[[t,t^{-1}]]/\C[[t]]$ is a differential module over the algebra $\C[t]$ with the derivation $\partial=d/dt$,   the notion of the  abstract  \emph{Martin integral} applies (see~\cite[sec. 2.3.1]{tv:git}) as the coset modulo the image of the derivation. 
 On $\C[[t,t^{-1}]]$, such an integral gives the residue at zero, up to a normalization factor, and on  $\C[[t]]$ it is zero. Hence on $M_t$   the induced integral will give   the rules~\eqref{eq.intdelta}.
%as desired. 

The difference with the usual delta-functions is in the equality $\delta(at)=a^{-1}\,\delta(t)$, without absolute value. Integral extends to several variables (with products of formal delta-functions with distinct arguments) and, basically, gives familiar formulas with determinants in place of absolute values of determinants. 

Indeed, consider first the single-variable case. Let $f(t)$ stand for a linear combination of formal $\delta(t)$ and its derivatives, $f(t)=\sum f_k\delta^{(k)}(t)$. Consider a simple change of variable $t=at'$ and $f(at')$ as a function of $t'$. Since for integration only the term $\delta(t')$ contributes and $\delta(t)=a^{-1}\delta(t')$, we obtain
\begin{equation*}
  \fint f(at')dt'=a^{-1}\fint f(t)dt\,.
\end{equation*}
If we write $t=t(t')$ for $t=at'$, so that $dt/dt'=a$, this can be re-written as
\begin{equation*}
  \fint f(t)dt = \fint f(t(t'))\frac{dt}{dt'}\,dt'\,.
\end{equation*}
For many-variables case, recall from Remark~\ref{rem.substdeltas} that the product of formal delta-functions produces $1/\det g$ under a non-degenerate linear substitution with a matrix $g$.  (Respectively, $1/\Ber g$ when there are also odd variables.) By the same argument, we obtain the formula
\begin{equation}
\label{eq.chvarfint}
  \fint Dt f(t)  = \fint Dt' \Ber\frac{\partial t}{\partial t'}\, f(t(t'))\,,
\end{equation}
where $t$ and $t'$ stand for   collections of even and odd variables. (In contrast with the usual Berezin integral, no extra signs related with orientation!)

\subsection{Formal Baranov-Schwarz transformation.}

Consider a formal generalized pseudo-differential form  or a \emph{formal pseudoform} for short. It is a linear combination of the products~\eqref{eq.fpdf}.  Denote it $\omega=\omega(dx)$. Here $dx$ stands for all $dx^a,d\xi^{\mu}$. Fix a dimension $r|s$, where $0\leq s\leq m$. (A priori $s$ does not have to be the same as in~\eqref{eq.fpdf}, but later we will see that it should be.) Mimicking Section~\ref{Sec:BSch}, consider the substitution $dx^A=dt^F\dot x_F^A$, where $\widetilde{dt^F}=\tilde{F}+1$. Define the \emph{formal Baranov-Schwarz transformation} (or \emph{transform}) $\lambda_{r|s}\omega$ as a function of $r$ even and $s$ odd variables collectively denoted $\dot x$, where
\begin{equation}\label{eq.fBS}
  (\lambda_{r|s}\omega)(\dot x):= \fint D(dt)\,\omega(dt^F\dot x_F)\,.
\end{equation}
The integral is the formal integral defined above, over $s$ even and $r$ variables.   We assume that the substitution into the formal pseudoform $\omega$ is well defined and maps it to a formal function of $dt^F$ of the type considered above, for the integral to make sense. %Below we will see how it works practically. 

\begin{theorem}
\label{thm.fBS}
If defined, the image $\lambda_{r|s}\omega$ of the formal Baranov-Schwarz transformation is a Voronov-Zorich $r|s$-form.
\end{theorem}

Note that, unlike before, the  formal Baranov-Schwarz transformation    takes values in ordinary Voronov-Zorich forms (not those of type II). 

The theorem will basically follow from the  two lemmas. Denote $L=\lambda_{r|s}\omega$.

\begin{lemma}
\label{lem.berfBS}
  $L=L(\dot x)$ satisfies condition~\eqref{eq.bercond}.
\end{lemma}
\begin{proof}
  Consider $L(g\dot x)$, where $g\in GL(r|s)$. We have
  \begin{multline*}
    L(g\dot x)=  \fint D(dt)\,\omega(\underbrace{dt^Fg_F^H}_{d{t'}^H}\dot x_H)=
    \fint D(dt'\cdot g^{-1})\,\omega(d{t'}^H\dot x_H)=\\
   \fint D(dt')\Ber((g^{-1})^{\Pi})\,\omega(d{t'}^H\dot x_H)= 
   \fint D(dt')\Ber g\,\omega(d{t'}^H\dot x_H)= \\
   \Ber g\cdot\fint D(dt')\,\omega(d{t'}^H\dot x_H)= \Ber g\cdot L(\dot x)\,,
  \end{multline*}
  where we used the change of variables formula for formal integral~\eqref{eq.chvarfint}. Here the parity-reversed matrix $(g^{-1})^{\Pi}$ arises because of the shift of parity for the indices of the variables $dt$ and $dt'$. Recall the property of   Berezinian $\Ber h^{\Pi}=(\Ber h)^{-1}$.  
\end{proof}

\begin{lemma}
\label{lem.fundfBS}
  $L=L(\dot x)$ satisfies condition~\eqref{eq.fund}.
\end{lemma}

This lemma can be proved either directly or by using the geometric meaning of condition~\eqref{eq.fund}. Namely, if $L=\lambda_{r|s}\omega$ is seen as a Lagrangian defining an action $S$ for $r|s$-paths  in a supermanifold in question (here we  consider $\omega(x,dx)$), then~\eqref{eq.fund} is equivalent to the condition that the variation $\delta S$ does not depend on the second derivatives~\cite{tv:pdf}.  

\begin{proof}[Proof of Lemma~\ref{lem.fundfBS}]
The variational derivative of an action $S$ with a Lagrangian $L=L(x,\dot x)$ is\,:
\begin{equation*}
  \frac{\delta S}{\delta x^A}= \frac{\partial L}{\partial x^A}- (-1)^{\tilde A\tilde F} \frac{\partial}{\partial t^F}\left(\frac{\partial L}{\partial \dot x_F^A}\right)
\end{equation*}
(see \cite{tv:git}). Hence, in our case, it will be (up to a common sign, which is not relevant),
\begin{multline*}
  \fint D(dt) \left[\frac{\partial \omega}{\partial x^A}-
  (-1)^{\tilde A\tilde F} \frac{\partial}{\partial t^F}\left(dt^F
  \frac{\partial \omega}{\partial dx^A}\right)(-1)^{(\tilde A+\tilde  F)(\tilde F+1)}
  \right]=\\
  \fint D(dt) \left[\frac{\partial \omega}{\partial x^A}-
  (-1)^{\tilde A} \frac{\partial}{\partial t^F}\left(dt^F
  \frac{\partial \omega}{\partial dx^A}\right)
  \right]=\\
   \fint D(dt) \left[\frac{\partial \omega}{\partial x^A}-
  \underline{(-1)^{\tilde A}
 dt^F} \left(\dot x_F^B\,
  \frac{\partial^2 \omega}{\partial x^B\partial dx^A}
  + \underline{\ddot x_{FG}^B\, dt^G\,
  \frac{\partial^2 \omega}{\partial dx^B\partial dx^A}(-1)^{(\tilde B +\tilde G)(\tilde G+1)}}
  \right)
  \right]\,.
\end{multline*}
The underlined terms are the only ones possibly containing the second derivatives of $x^A$. After simplifying the signs, they will be:
\begin{equation*}
  (-1)^{\tilde F(\tilde G+1)}dt^Fdt^G\,\ddot x_{FG}^B\,(-1)^{\tilde A}\omega_{BA}\,,
\end{equation*}
where $\omega_{BA}$ are the partial derivatives of $\omega$ in $dx^B$ and $dx^A$. So, effectively, we need to analyze $(-1)^{\tilde F(\tilde G+1)}dt^Fdt^G\,\ddot x_{FG}^B$. We have:
\begin{multline*}
  (-1)^{\tilde F(\tilde G+1)}dt^Fdt^G\,\ddot x_{FG}^B=
  (-1)^{\tilde G(\tilde F+1)}dt^Gdt^F\,\ddot x_{GF}^B=\\
  (-1)^{\tilde G(\tilde F+1)}(-1)^{(\tilde F+1)(\tilde G+1)}dt^Fdt^G\,(-1)^{\tilde F\tilde G}\ddot x_{FG}^B=
  -  (-1)^{\tilde F(\tilde G+1)}dt^Fdt^G\,\ddot x_{FG}^B=0\,,
\end{multline*}
as desired.
%(basically by the opposite commutativity rules for   differentials and partial derivatives).
\end{proof}

(Morally, $\delta S$ contains no second derivatives    because of the de Rham differential acting on formal pseudoforms and   a version of Stokes' formula,  underpinning the computation above.There is no difference here with the case of usual pseudoforms considered  in~\cite{tv:pdf}.) 

Together, Lemma~\ref{lem.berfBS} and Lemma~\ref{lem.fundfBS} give Theorem~\ref{thm.fBS}. 

So far, in our considerations we simply assumed that a formal Baranov-Schwarz transformation of a given formal pseudoform exists. Now we shall look at this closer. First of all, if we want the substitution $dx^A=dt^F\dot x_F^A$ to produce a function of the desired form in variables $dt^F$, so that formal integration to make sense, we must assume that $s$ in $\lambda_{r|s}$ is the same $s$ as for $\omega=\omega(dx)$, a linear combination of~\eqref{eq.fpdf},  to which the transformation will be applied. %So, we take $\omega=\omega(dx)$ which is a linear combination of~\eqref{eq.fpdf} and apply to it  $\lambda_{r|s}$ with exactly  the same $s$. 

For a fixed $s$, the space of formal pseudoforms is graded by $N\in\Z$, where
\begin{equation*}
  N=k_1+\ldots +k_n+i_1+\ldots + i_{m-s}-j_1-\ldots-j_s-s\,,
\end{equation*}
$k_1,\ldots,k_n\in\{0,1\}$ and $i_1,\ldots, i_{m-s}, j_1,\ldots, j_{s}\in \Z_{\geq 0}$. % and $j_1,\ldots, j_{s}\in \Z_{\geq 0}$. 
The number $N$ is the weight of $\omega(dx)$ with respect to scaling transformations $dx\mapsto a\cdot dx$.

The following statement is immediate and applies both to usual and formal pseudoforms.

\begin{theorem}
  Suppose a pseudo-differential form $\omega(dx)$   is homogeneous of weight $N$ under scaling $dx\mapsto a\,dx$. Then its (formal or usual) Baranov-Schwarz transformation $\lambda_{r|s}\omega$ is either zero or undefined if $N\neq r-s$.
\end{theorem}
\begin{proof}
  Suppose $L=\lambda_{r|s}\omega$. By the definition of $\lambda_{r|s}$, $L(a\dot x)=\lambda_{r|s}[\omega(a\,dx)]$. Hence, if $\omega$ is homogeneous of weight $N$ in $dx$, then $L$ has to be homogeneous of the same weight $N$ in $\dot x$. But by the definition of $r|s$-forms, $L$ is homogeneous of weight $r-s$. So, $N=r-s$ unless $L\equiv 0$\,.
\end{proof}

This statement is analogous to the fact that, for ordinary manifolds, where pseudoforms are inhomogeneous differential forms, the Baranov-Schwarz transformation is the projector onto the subspace of forms of a particular degree\footnote{Where in the source forms are treated as functions of differentials and in the target as functions of tangent vectors.}.  

Therefore we obtain the following relation between the weight $N$ for a formal pseudoform (with a fixed $s$)  and the degree $r$ of the $r|s$-form corresponding to it under the transformation $\lambda_{r|s}$\,:
\begin{equation}\label{eq.Nr}
  N=r-s\,.
\end{equation}
In algebraic setting, we conclude that the spaces $S^{N,s}(\Pi V)$ introduced in Section~\ref{sec.geom} heuristically correspond to the super exterior powers $\Lambda^{N+s|s}(V)$, up to $\Pi^{N}=\Pi^{r-s}$. 

\begin{remark}
In the definition of $r|s$-forms  that we have used, we have $r\geq 0$, as the number of even vectors in the argument. (The same for    the corresponding notion of $r|s$-vectors, where the arguments are $r$ even and $s$ covectors, respectively.) However, in~\eqref{eq.Nr}, if $N$ is negative and $|N|>s$, then $r<0$. Therefore, to get a full theory, one must make use of $r|s$-forms with negative $r$ introduced in~\cite{tv:dual}. Geometrically they correspond to replacing ``paths'' by a dual notion of ``copaths'' (which are systems of equations). One can introduce  a version of formal Baranov-Schwarz transformation for such a dual setting. We will consider such an extended theory elsewhere. 
\end{remark}

\begin{theorem}
\label{thm.dBSmono}
  For a fixed $s$ and $N=r-s$, the formal Baranov-Schwarz transformation is well defined and injective. 
\end{theorem}

We omit the proof of this theorem, which will be given in our subsequent paper. In view of Theorem~\ref{thm.dBSmono}, it is tempting to identify, up to a parity reversion, the space of $r|s$-forms  with the image of formal $\lambda_{r|s}$. (In algebraic setting, the space $\Lambda^{r|s}(V)$ with $\Pi^{r-s}S^{r-s,s}(\Pi V)$.)

If we assume~\eqref{eq.Nr} and use $r$ instead of $N$ as a $\Z$-grading, then, for a fixed $s$, the general appearance of formal pseudoforms of degree $r$ will be:
\begin{equation}\label{eq.fpdfdegr}
  \omega=  
	\sum_{\substack{k+i-j=r\\
k_1+\ldots+k_n=k\\
i_1+\ldots+i_{m-s}=i\\
			j_{1}+\ldots+j_s=j}} \!\!
  (dx^1)^{k_1}\,\ldots\, (dx^n)^{k_n}\,
	(d\xi^1)^{i_1}\ldots (d\xi^{m-s})^{i_{{m-s}}}\,
	\delta^{(j_1)}(d\xi^{m-s+1})\,\ldots \,\delta^{(j_s)}(d\xi^m)\;
\omega_{\mathrm{K}\mathrm{I}|\mathrm{J}}\,,
%\omega_{k_1,\ldots,k_ni_1,\ldots,i_{m-s}|j_{1},\ldots,j_s}
\end{equation}
where $\mathrm{K}$, $\mathrm{I}$, $\mathrm{J}$ are multi-indices, %$\mathrm{K}=(k_1,\ldots,k_n)$, $\mathrm{I}=(i_1,\ldots,i_{m-s})$ and $\mathrm{J}=(j_{1},\ldots,j_s)$. 
 $k_1,\ldots,k_n\in \{0,1\}$  and $i_1,\ldots,i_{m-s},j_{1},\ldots,j_s\in \Z_{\geq 0}$.

\subsection{Examples.}

To show how formal Baranov-Schwarz transformation works, we can consider the simplest non-trivial case of dimension $n|m=1|2$. Then $s$ can be $0$, $1$ or $2$, and $\omega$ will be, respectively:
\begin{equation}\label{eq.oms0}
  \omega=\sum\limits_{k+i_1+i_2=r} (dx)^k\,(d\xi^1)^{i_1}\,(d\xi^2)^{i_2}\;\omega_{ki_1i_2}%=
  %\sum\limits_{i_1+i_2=r} (d\xi^1)^{i_1}(d\xi^2)^{i_2}\;\omega_{0i_1i_2}+
  %dx\sum\limits_{i_1+i_2=r-1} (d\xi^1)^{i_1}(d\xi^2)^{i_2}\;\omega_{1i_1i_2}
\end{equation}
for $s=0$; or
\begin{equation}\label{eq.oms1}
  \omega=\sum\limits_{k+i-j=r} (dx)^k\,(d\xi^1)^{i}\,\delta^{(j)}(d\xi^2)\;\omega_{ki|j}
\end{equation}
for $s=1$; or
\begin{equation}\label{eq.oms2}
  \omega=\sum\limits_{k-j_1-i_2=r} (dx)^k\,\delta^{(j_1)}(d\xi^1)\,\delta^{(j_2)}(d\xi^2)\;\omega_{k|j_1j_2}
\end{equation}
for $s=2$. In case~\eqref{eq.oms0}, which is just the ordinary differential forms, the space of such forms is finite-dimensional for a given $r$, $r\geq 0$ and $r$ is not bounded above.   In case~\eqref{eq.oms2}, which would correspond to integral forms, the space is again finite-dimensional for a given $r$, where $r\leq 1$ and  not bounded below.   For the intermediate case~\eqref{eq.oms1}, we can re-write
\begin{multline}\label{eq.oms1det}
  \omega=\sum\limits_{k+i-j=r} (dx)^k\,(d\xi^1)^{i}\,\delta^{(j)}(d\xi^2)\;\omega_{ki|j}= \\
  \sum\limits_{i-j=r}  (d\xi^1)^{i}\,\delta^{(j)}(d\xi^2)\;\omega_{0i|j}+
  dx\sum\limits_{i-j=r-1} (d\xi^1)^{i}\,\delta^{(j)}(d\xi^2)\;\omega_{1i|j}\,.
\end{multline}
Here $r$ can be any integer. If $r\geq 0$, we can further express $\omega$ as  
\begin{equation}\label{eq.oms1fdet}
  \omega=
  \sum\limits_{j\geq 0}  (d\xi^1)^{j+r}\,\delta^{(j)}(d\xi^2)\;\omega_{0|j}+
  dx\sum\limits_{j\geq 0} (d\xi^1)^{j+r-1}\,\delta^{(j)}(d\xi^2)\;\omega_{1|j}\,.
\end{equation}
If $r< 0$, we can similarly write 
\begin{equation}\label{eq.oms1fdetneg}
  \omega=
  \sum\limits_{i\geq 0}  (d\xi^1)^{i}\,\delta^{(i-r)}(d\xi^2)\;\omega_{0i}+
  dx\sum\limits_{i\geq 0} (d\xi^1)^{i}\,\delta^{(i-r+1)}(d\xi^2)\;\omega_{1i}\,.
\end{equation}
In both cases, for a given $r\in \Z$, the space will be infinite-dimensional. 

%For the purpose of formal Baranov-Schwarz transformation $\lambda_{r|s}$, 

We now consider $r\geq 0$.  

We skip $\lambda_{r|0}$ applied to~\eqref{eq.oms0} because it would give just the familiar isomorphism of polynomials in $dx,d\xi^1,d\xi^2$ of degree $r$ with multilinear antisymmetric functions of $r$ vectors and consider $s=2$ and $s=1$. 

Case $s=2$ is close to what we considered in Section~\ref{Sec:BSch}. For $r\geq 0$, $r=0$  or $r=1$.

\begin{example}[$\lambda_{0|2}$ in dimension $1|2$]
We have
\begin{equation}\label{eq.02in12}
  \omega=\delta(d\xi^1)\delta(d\xi^2)\,\omega_{0|00}+
  dx\Bigl(\delta(d\xi^1)\delta'(d\xi^2)\,\omega_{1|01}+ \delta'(d\xi^1)\delta(d\xi^2)\,\omega_{1|10}\Bigr)\,.
\end{equation}
Substitution: $dx=d\tau^1\dot x_{\hat 1}+ d\tau^2\dot x_{\hat 2}$, $d\xi^1=d\tau^1\dot \xi^1_{\hat 1}+ d\tau^2\dot \xi^1_{\hat 2}$ and 
$d\xi^2=d\tau^1\dot \xi^2_{\hat 1}+ d\tau^2\dot \xi^2_{\hat 2}$. To evaluate $\lambda_{0|2}\omega$, we will need to take the formal integral $\fint D(d\tau^1,d\tau^2)$. Instead of the direct computation, we note that the matrix $(\dot\xi^{\mu}_{\hat \nu})$ is invertible, hence it is possible to consider $d\xi^1,d\xi^2$ as variables of integration. Then $dx=d\xi^1\,\alpha+ d\xi^2\,\beta$, where $\left[\begin{smallmatrix} \alpha \\ \beta \end{smallmatrix} \right]=(\dot\xi^{\mu}_{\hat \nu})^{-1}\left[\begin{smallmatrix} \dot x_{\hat 1} \\ \dot x_{\hat 2}  \end{smallmatrix}\right]$. We obtain:
\begin{equation*}
  dx\, \delta(d\xi^1)\delta'(d\xi^2)= (d\xi^1\,\alpha+ d\xi^2\,\beta)\,\delta(d\xi^1)\delta'(d\xi^2)=
  d\xi^2\,\beta\,\delta(d\xi^1)\delta'(d\xi^2)= 
  -\delta(d\xi^1)\delta(d\xi^2)\,\beta\,.
\end{equation*}
Similarly, $dx\, \delta'(d\xi^1)\delta(d\xi^2)= -\delta(d\xi^1)\delta(d\xi^2)\,\alpha$\,. Thus we have after the substitution
\begin{equation*}
  \omega=\delta(d\xi^1)\delta(d\xi^2)\,\bigl(\omega_{0|00}-\beta\omega_{1|01}-\alpha\omega_{1|10}\bigr)=
  \delta(d\tau^1)\delta(d\tau^2)\,\frac{1}{\det(\dot\xi^{\mu}_{\hat \nu})}\,\bigl(\omega_{0|00}-\beta\omega_{1|01}-\alpha\omega_{1|01}\bigr)\,.
\end{equation*}
Recalling the definition of $\alpha,\beta$ and the formula for Berezinian, we can express it as 
 \begin{equation*}
   \omega= \delta(d\tau^1)\delta(d\tau^2)\,\Ber\!
   \left(
\begin{array}{c|cc}
	\omega_{0|00} & -\omega_{1|10} & -\omega_{1|01} \\
	\hline
	\dot x_{\hat 1} & \dot \xi^1_{\hat 1} & \dot \xi^2_{\hat 1} \\
	\dot x_{\hat 2} & \dot \xi^1_{\hat 2} & \dot \xi^2_{\hat 2} 
\end{array}
\right)\,.
 \end{equation*}
Hence, finally, after integration, we obtain:
\begin{equation}\label{eq.BS02in12}
  (\lambda_{0|2}\omega)\begin{pmatrix}
                         \dot x_{\hat 1} & \dot \xi^1_{\hat 1} & \dot \xi^2_{\hat 1} \\
	\dot x_{\hat 2} & \dot \xi^1_{\hat 2} & \dot \xi^2_{\hat 2}
                       \end{pmatrix}=
                       \Ber\!
   \left(
\begin{array}{c|cc}
	\omega_{0|00} & -\omega_{1|10} & -\omega_{1|01} \\
	\hline
	\dot x_{\hat 1} & \dot \xi^1_{\hat 1} & \dot \xi^2_{\hat 1} \\
	\dot x_{\hat 2} & \dot \xi^1_{\hat 2} & \dot \xi^2_{\hat 2} 
\end{array}
\right)\,.
\end{equation}
(Compare with Example~\ref{ex.BS01in11} in Section~\ref{Sec:BSch}.)  
\end{example}

\begin{example}[$\lambda_{1|2}$ in dimension $1|2$]
We have
\begin{equation}\label{eq.12in12}
  \omega=dx\,\delta(d\xi^1)\delta(d\xi^2)\,\omega_{1|00}
\end{equation}
and
\begin{equation}\label{eq.BS12in12}
  (\lambda_{1|2}\omega)\left(
\begin{array}{c|cc}
	\dot x_{1} & \dot \xi^1_{1} & \dot \xi^2_{1} \\
	\hline
	\dot x_{\hat 1} & \dot \xi^1_{\hat 1} & \dot \xi^2_{\hat 1} \\
	\dot x_{\hat 2} & \dot \xi^1_{\hat 2} & \dot \xi^2_{\hat 2} 
\end{array}
\right)=
                       \Ber\!
   \left(
\begin{array}{c|cc}
	\dot x_{1} & \dot \xi^1_{1} & \dot \xi^2_{1} \\
	\hline
	\dot x_{\hat 1} & \dot \xi^1_{\hat 1} & \dot \xi^2_{\hat 1} \\
	\dot x_{\hat 2} & \dot \xi^1_{\hat 2} & \dot \xi^2_{\hat 2} 
\end{array}
\right)\omega_{1|00}\,.
\end{equation}
\end{example}

Consider now the  case $s=1$. Let $\omega$ be given by~\eqref{eq.oms1fdet}. We will take $r=1$ and $r=0$ for simplicity.  
 
\begin{example}[$\lambda_{1|1}$ in dimension $1|2$]
Here
\begin{equation}\label{eq.om11}
  \omega=
  \sum\limits_{j\geq 0}  (d\xi^1)^{j+1}\,\delta^{(j)}(d\xi^2)\;\omega_{0|j}+
  dx\sum\limits_{j\geq 0} (d\xi^1)^{j}\,\delta^{(j)}(d\xi^2)\;\omega_{1|j}\,.
\end{equation}
Substitution: $dx=dt\cdot \dot x_1+d\tau\cdot \dot x_{\hat 1}$, $d\xi^1=dt\cdot \dot \xi^1_1+d\tau\cdot \dot \xi^1_{\hat 1}$, and $d\xi^2=dt\cdot \dot \xi^2_1+d\tau\cdot \dot \xi^2_{\hat 1}$. As above, we assume that it is non-degenerate, which here means $\dot \xi^2_{\hat 1}$ is invertible. So, we take $dt$ and $d\xi^2$ as variables instead of $dt$ and $d\tau$. We obtain $dx=dt\cdot a +d\xi^2\cdot \beta$ and $d\xi^1=dt\cdot \gamma+ d\xi^2\cdot c$, where $a=\dot x_1-(\dot\xi^2_{\hat 1})^{-1}\dot\xi^2_1\dot x_{\hat 1}$, $\beta=(\dot\xi^2_{\hat 1})^{-1}\dot x_{\hat 1}$, $\gamma=\dot\xi^1_1-(\dot\xi^2_{\hat 1})^{-1}\dot\xi^2_1\dot\xi^1_{\hat 1}$ and $c=(\dot\xi^2_{\hat 1})^{-1}\dot\xi^1_{\hat 1}$. In this notation, we obtain that upon substitution
\begin{equation*}
  (d\xi^1)^{j+1}\,\delta^{(j)}(d\xi^2)= \delta (d\xi^2) dt\cdot (-1)^j(j+1)! c^j \gamma
\end{equation*}
and
\begin{equation*}
  dx\,(d\xi^1)^{j}\,\delta^{(j)}(d\xi^2)= \delta (d\xi^2) dt\cdot (-1)^jj!(-ac^j+j\beta\gamma)
\end{equation*}
(the calculation uses the identity $t^n\delta^{(n)}(t)=(-1)^nn!\delta(t)$ for formal delta-functions). Notice that $\delta (d\xi^2) dt=(\dot\xi^2_{\hat 1})^{-1}\cdot \delta(d\tau)dt$ (which we need for integration). Skipping intermediate calculations, we can give the final answer as:
\begin{multline}\label{eq.BS11in12}
  (\lambda_{1|1}\omega)\left(
\begin{array}{c|cc}
	\dot x_{1} & \dot \xi^1_{1} & \dot \xi^2_{1} \\
	\hline
	\dot x_{\hat 1} & \dot \xi^1_{\hat 1} & \dot \xi^2_{\hat 1} 
\end{array}
\right)=
\sum\limits_{j\geq 0}(-1)^j(j+1)!\,(\dot\xi^2_{\hat 1})^{-j}(\dot\xi^1_{\hat 1})^j\cdot 
\Ber\!\left(
\begin{array}{c|c}
	 \dot \xi^1_{1} & \dot \xi^2_{1} \\
	\hline
	\dot \xi^1_{\hat 1} & \dot \xi^2_{\hat 1} 
\end{array}
\right)\omega_{0|j}+\\
\sum\limits_{j\geq 0}(-1)^j j!\,
\left(   
-(\dot\xi^2_{\hat 1})^{-j}(\dot\xi^1_{\hat 1})^j
\cdot 
\Ber\!\left(
\begin{array}{c|c}
	 \dot x_{1} & \dot \xi^2_{1} \\
	\hline
	\dot  x_{\hat 1} & \dot \xi^2_{\hat 1} 
\end{array}
\right) 
+
j (\dot\xi^2_{\hat 1})^{-1} \dot x_{\hat 1}  \cdot
\Ber\!\left(
\begin{array}{c|c}
	 \dot \xi^1_{1} & \dot \xi^2_{1} \\
	\hline
	\dot  \xi^1_{\hat 1} & \dot \xi^2_{\hat 1} 
\end{array}
\right) 
\right)\omega_{1|j}\,.
\end{multline}
Note that there are two Berezinian minors of the matrix in the argument here, one of them being what we call ``wrong Berezinian''~\cite{shemya:voronov2019:superplucker},\cite{shemyakova:2023}, which is a Berezinian with an odd column   or row vector inserted into an even position (and whose value is therefore odd). 
\end{example}

(It would be interesting to compare this formula~\eqref{eq.BS11in12} with the description in Section~\ref{Sec:11}.)
 
\begin{example}[$\lambda_{0|1}$ in dimension $1|2$]
Here
\begin{equation}\label{eq.om01}
  \omega=
  \sum\limits_{j\geq 0}  (d\xi^1)^{j}\,\delta^{(j)}(d\xi^2)\;\omega_{0|j}+
  dx\sum\limits_{j\geq 1} (d\xi^1)^{j-1}\,\delta^{(j)}(d\xi^2)\;\omega_{1|j}\,.
\end{equation}
  Acting as before, we consider the substitution $dx = d\tau \dot x_{\hat 1}$, $d\xi^1=d\tau \dot \xi^1_{\hat 1}$, and $d\xi^2=d\tau \dot \xi^2_{\hat 1}$. We need to assume   $\dot\xi^2_{\hat 1}$ to be invertible in order for the formal delta-function to be defined. By using the same identities as in the previous example and taking the formal integral over $d\tau$ (where $\delta(d\xi^2_{\hat 1})=(\dot \xi^2_{\hat 1})^{-1}\delta(d\tau)$), we will arrive at
\begin{equation}\label{eq.BS01in12}
   (\lambda_{1|1}\omega)%\left
   (\dot x_{\hat 1}\,|\, \dot \xi^1_{\hat 1},   \dot \xi^2_{\hat 1} %\right
   )=
\sum\limits_{j\geq 0}(-1)^jj!\,(\dot \xi^2_{\hat 1})^{-j-1}(\dot \xi^1_{\hat 1})^{j}\,\omega_{0|j}-
\dot x_{\hat 1}\sum\limits_{j\geq 1}(-1)^jj!\,(\dot \xi^2_{\hat 1})^{-j-1}(\dot \xi^1_{\hat 1})^{j-1}\,\omega_{1|j}\,.
\end{equation}
This is essentially a general form of a homogeneous function of degree $-1$ of an odd vector $(\dot x_{\hat 1}\,|\, \dot \xi^1_{\hat 1},   \dot \xi^2_{\hat 1})$, provided $\dot\xi^2_{\hat 1}$ is invertible. (Recall that the homogeneity condition is all that remains from the definition of $r|s$-form for $r|s=0|1$.)
\end{example}
 
A common feature  of the above examples is that the $r|s$-form obtained as the formal Baranov-Schwarz transformation is a rational function  with its pole determined by the splitting of the odd coordinates $\xi^{\mu}$ into the groups $\xi^1,\ldots,\xi^{m-s}$ and $\xi^{m-s+1},\ldots,\xi^{m}$. % that we have adopted. 

Finally,  we want to consider an example of a more general type of formal pseudoforms. % than   considered before. 

\begin{example}\label{ex.last}
Take a collection  of $r$ odd and $s$ even $1$-forms (which is equivalent to $r$ even and $s$ odd covectors). Denote them  $\alpha^1,\ldots,\alpha^r$, $\alpha^{\hat 1},\ldots,\alpha^{\hat s}$. We can write them as $\alpha^F$,  where $F=1,\ldots,r,\hat 1, \ldots, \hat s$ and $\widetilde{\alpha^F}=\tilde F+1$. If 
\begin{equation*}
  \alpha^F=dx^A\alpha_A^F\,, 
\end{equation*}
then the matrix $(\alpha_A^F)$ will be even, i.e. $\widetilde{\alpha_A^F}=\tilde{A}+\tilde{F}$.  Define a generalized pseudoform as
\begin{subequations}\label{eq.nwpalphas}
\begin{equation}
  \omega:= \alpha^1 \ldots \alpha^r \,\delta(\alpha^{\hat 1})\ldots \delta(\alpha^{\hat s}) \equiv \delta(\alpha)\,.
  \tag{\ref{eq.nwpalphas}}
\end{equation}
A particular special case: $\alpha^1=dx^1, \ldots, \alpha^r=dx^r$, $\alpha^{\hat 1}=d\xi^{m-s+1}, \ldots, \alpha^{\hat s}=d\xi^{m}$. Then 
\begin{equation} \label{eq.nwpalphasA}
  \omega = dx^1\ldots dx^r\,\delta(d\xi^{m-s+1})\ldots \delta(d\xi^{m})\,.
\end{equation}
Another special case: $\alpha^1=dx^{a_1}, \ldots, \alpha^r=dx^{a_r}$, $a_1<\ldots a_r$, $\alpha^{\hat 1}=d\xi^{\mu_1}, \ldots, \alpha^{\hat s}=d\xi^{\mu_s}$, $\mu_1<\ldots \mu_s$. Then
\begin{equation} \label{eq.nwpalphasB}
  \omega = dx^{a_1}  \ldots dx^{a_r}\,\delta(d\xi^{\mu_1})\ldots \delta(d\xi^{\mu_s})\,.
\end{equation}
 \end{subequations} 
If we substitute $dx^A=dt^G\dot x^A_G$ into $\omega$, then for each $\alpha^F$ we will get $\alpha^F=dt^G\dot x^A_G\alpha_A^F$, hence
\begin{subequations}\label{eq.BSpalphas}
\begin{equation}
  (\lambda_{r|s}\omega)(\dot x)=\fint D(dt)\,\delta (dt\,\dot x^A \alpha_A\,)= \Ber(\dot x^A_G \alpha_A^F)= \Ber\left(\langle \dot x_G,\alpha ^F\rangle\right)\,.
  \tag{\ref{eq.BSpalphas}}
\end{equation}
(The Berezinian appears in the numerator because of the parity reversion for the indices.) The $r|s$-form $\lambda_{r|s}\omega$ is a ``non-linear analogue of the wedge product'' for the $1$-forms $\alpha^F$ (compare in~\cite{shemya:voronov2019:superplucker}).  For the special cases above, we would obtain minors of the matrix $\dot x_G^F$:
\begin{equation}\label{eq.minorA}
  (\lambda_{r|s}\omega)(\dot x)=\Ber\!
   \left(
\begin{array}{ccc|ccc}
	\dot x_{1}^1& ...& \dot x_1^r& \dot \xi_1^{m-s+1} & ... & \dot \xi_1^m \\
...& ...& ...& ... & ... & ... \\
\dot x_{r}^1& ...& \dot x_r^r& \dot \xi_r^{m-s+1} & ... & \dot \xi_r^m \\
	\hline
	\dot x_{\hat 1}^1& ...& \dot x_{\hat 1}^r& \dot \xi_{\hat 1}^{m-s+1} & ... & \dot \xi_{\hat 1}^m \\
...& ...& ...& ... & ... & ... \\
\dot x_{\hat s}^1& ...& \dot x_{\hat s}^r& \dot \xi_{\hat s}^{m-s+1} & ... & \dot \xi_{\hat s}^m \\
\end{array}
\right)
\end{equation}
or
\begin{equation}\label{eq.minorB}
  (\lambda_{r|s}\omega)(\dot x)=
  \Ber\!
   \left(
\begin{array}{ccc|ccc}
	\dot x_{1}^{a_1}& ...& \dot x_1^{a_r}& \dot \xi_1^{\mu_1} & ... & \dot \xi_1^{\mu_s} \\
...& ...& ...& ... & ... & ... \\
\dot x_{r}^{a_1}& ...& \dot x_r^{a_r}& \dot \xi_r^{\mu_1} & ... & \dot \xi_r^{\mu_s} \\
	\hline
	\dot x_{\hat 1}^{a_1}& ...& \dot x_{\hat 1}^{a_r}& \dot \xi_{\hat 1}^{\mu_1} & ... & \dot \xi_{\hat 1}^{\mu_s} \\
...& ...& ...& ... & ... & ... \\
\dot x_{\hat s}^{a_1}& ...& \dot x_{\hat s}^{a_r}& \dot \xi_{\hat s}^{\mu_1} & ... & \dot \xi_{\hat s}^{\mu_s} \\
\end{array}
\right)
\end{equation}
 \end{subequations} 
\end{example} 

The above considerations make  sense if the odd-odd block of $(\alpha_A^F)$ has rank $s$, i.e. the $1$-forms $\alpha^{\hat 1}, \ldots, \alpha^{\hat s}$ are independent and span an $s$-dimensional subspace in $\Pi T_x^*M$.  
The   formal pseudoform given by~\eqref{eq.nwpalphas} is concentrated on its annihilator. For example, in the case~\eqref{eq.nwpalphasA}, it is given by the equations $d\xi^{m-s+1}=0, \ldots, d\xi^m=0$. Respectively, the  
 $r|s$-form that we obtain has  a pole whenever any of the odd arguments $\dot x_{\hat 1}, \ldots, \dot x_{\hat s}$ is in this subspace. 

If we fix $\omega$ of the form~\eqref{eq.nwpalphas} and act on it by multiplication by various $1$-forms and by the ``inner products'' $i_u=(-1)^{\tilde u}u^A\partial/\partial dx^A$ with tangent vectors $u=u^A\partial_A$, we will obtain a Fock space similar to that we considered above (and coinciding with it for $\omega$ given by~\eqref{eq.nwpalphasA}). The ``vacuum vector'' there will be formal pseudoform $\delta(\alpha^{\hat 1})\ldots \delta(\alpha^{\hat s})$.

\section{Conclusion.}

We   investigated Voronov-Zorich forms and their algebraic version, the super exterior powers $\Lambda^{r|s}(V)$ of $n|m$-dimensional vector superspace. The main problem   is lack of an explicit description except for the cases $s=0$ and $s=m$. 
For $\Lambda^{1|1}(V)$, we established an isomorphism with closed differential $1$-forms on a projective superspace  $\mathbb{P}^{m-1|n}$. Then we studied the Laurent expansions   of the function $R_A(z)=\Ber(E+zA)$, where $A$ is  an even linear operator. It was found previously  (Schmitt and Khudaverdian-Voronov) that the Taylor expansion  of  $R_A(z)$  at $z=0$ encodes  ordinary exterior powers $\Lambda^k(V)$ (as the supertraces of the action of $A$), while the expansion near $z=\infty$  encodes similarly the algebraic analogue of the Bernstein-Leites integral forms. We   found here  that a series  of  vector spaces  constructed explicitly in terms of   bases correspond  in the same way to the Laurent expansions in all annular regions  between the poles. The number of such   regions together with the neighborhoods of zero and infinity is $m$,  and the constructed spaces  model all the sought-for different types of the exterior powers $\Lambda^{r|s}(V)$ for $s=0,\ldots,m$. More precisely, they should be viewed as ``generalized symmetric powers'' $S^{N,s}(\Pi V)$.  As a tool, we used the language of formal delta-functions that appeared before in physics in the works of Belopolsky and Witten. In differential-geometric framework, we showed that the corresponding formal pseudoforms are   related with Voronov-Zorich $r|s$-forms  by  a formal version of the Baranov-Schwarz integral transformation. Next steps that we will consider elsewhere should  include the study of coordinate transformation properties  and  incorporating into the  picture of  $r|s$-forms for   $r<0$ (introduced in~\cite{tv:dual}). Remarks after Example~\ref{ex.last} indicate at links with geometry of the Grassmannian, integral geometry in the sense of Gelfand-Gindikin-Graev (first discovered in~\cite{tv:pdf, tv:cohom-1988}) and Gelfand's theory of general hypergeometric functions (on it, see~\cite{gelfand-graev-retakh1992}), which all need to be further investigated. 

%\bibliographystyle{plain}
%\bibliography{../../../../general.bib}{}
%\bibliography{general.bib, geometry.bib}{}
%\end{document}

\end{document}